\documentclass[a4paper]{article}

\usepackage{amssymb,amsmath,amsfonts,latexsym}
\usepackage[latin1]{inputenc}
\usepackage[T1]{fontenc}
\usepackage{graphicx}
\usepackage{eurosym, multicol}
\usepackage[english]{babel}
\usepackage{amsthm}
\usepackage{tikz}
  \usepackage[all]{xy}
  \usepackage{bm}
  \usetikzlibrary{matrix}
\usepackage{hyperref}
\usepackage{enumerate}

\newcommand{\N}{\mathbb{N}}
\newcommand{\nin}{\not\in}
\newcommand{\nni}{\not\ni}
\newcommand {\Sub}{Sub}

\DeclareMathOperator{\Of}{Clop}
\newcommand{\R}{\mathrel{R}}
\DeclareMathOperator{\Ult}{Ult}
\DeclareMathOperator{\cl}{cl}
\DeclareMathOperator{\Con}{Con}
\DeclareMathOperator{\Var}{Var}
\DeclareMathOperator{\id}{id}
\DeclareMathOperator{\Log}{Log}
\DeclareMathOperator{\At}{At}
\newcommand{\K}{\mathcal{K}}
\DeclareMathOperator{\mo}{mod}

\newcommand{\recprec}{\mathrel{\prec \!\!\!\!\!\!^\smallsmile}}

\usepackage{tikz}
  \usepackage[all]{xy}
  \usepackage{bm}
  \usetikzlibrary{matrix}

\theoremstyle{definition}
\newtheorem{definition}{Definition}[section]
\newtheorem{remark}[definition]{Remark}
\newtheorem{example}[definition]{Example}
\newtheorem{nota}[definition]{Notation}

\theoremstyle{plain}
\newtheorem{theorem}[definition]{Theorem}
\newtheorem{lemma}[definition]{Lemma}
\newtheorem{corollary}[definition]{Corollary}
\newtheorem{proposition}[definition]{Proposition}
\newtheorem{problem}[definition]{Problems}

\title{Subordination Algebras in Modal Logic}

  \author{Laurent De Rudder, Georges Hansoul and Valentine Stetenfeld}
 \date{}

\begin{document}
\maketitle

{\abstract{  Subordination algebras are Boolean algebras together with a binary relation subject to some conditions. They proved to be useful for a region based study of space (\cite{Proximity}) and for de Vries duality for compact Hausdorff spaces (\cite{deVries}). In both cases, the relation has a clear geometrical flavour. They have been studied by Celani (\cite{Quasimod}, \cite{qcong}) under the name "quasi-modal algebras" in a more universal algebraic and logic spirit: the binary relation is considered as a (multi-valuated) operator which obeys the same equations as standard modal operators.
  
  This paper goes one step further in this direction. In Section \ref{Section1}, we recall and use Celani's duality results to define two modal operators on the canonical extension of a subordination algebra. This enables, in Section \ref{Section2}, to define validity of modal (or bimodal) formulas on any subordination algebra and obtain completeness results. We show in Section \ref{Section_universal algebra} that Celani's notions of subalgebras and homomorphic image preserves validity, and examine the corresponding problem for products, to obtain a Birkhoff HSP-theorem. Finally, in Section \ref{Section3}, we extend to the new setting correspondence theory: we adapt Sahlqvist's theorem and undertake a study of the correspondence between subordination and modal (or even bimodal) formulas. } 
  
  }
  
  
  
\section*{Introduction}
Modal algebra is a powerful tool of investigation of normal modal logics, which is dual, in the sense of Stone duality, to the Kripke semantic.

Apparently far from the modal land, a contact algebra (see for instance \cite{ContactDimov}, \cite{Contact} or \cite{Koppelberg}) is a hybrid structure (both algebraic and relational) useful in the spatial reasoning in terms of regions rather than in terms of points. Moreover, they are a major ingredient of de Vries celebrated duality for compact Hausdorff spaces (\cite{deVries}).

Despite the very distinct domains of application, we have here very similar mathematical objects. Both are Boolean algebras with an extra structure: a modal operator in the first case, a binary relation (contact or proximity) in the second case. A the level of Stone duality, the similarity is even clearer: the dual of a modal algebra is a Stone space with a closed "accessibility" relation which is continuous in the sense that the inverse image of a clopen subset is clopen, while the dual of a contact algebra is a closed "adjency" relation which is reflexive and symmetric. 

There is a common denominator to these topological structures: Stone spaces with a closed relation. These are duals (see \cite{Quasimod}) of algebras that have been introduced and studied under several denominations: subordination algebras \cite{Sourab}, pre-contact or proximity algebras \cite{PrecontactDimov} and \cite{Proximity} and quasi-modal algebras \cite{Quasimod}. On these objects (and it does not matter to consider the algebraic, syntactical side, or its relational, semantical side) coexist three distinct language: the modal language (say the diamond $\lozenge$), the contact language (subordination $\prec$ or contact $\mathcal{C}$) and the relational or Kripke language (accessibility $R$). We propose in this paper to focus in the relations between the subordination world and the modal one.

Celani's point of view in \cite{Quasimod} seems to be a good start in this direction: subordination algebras are presented as a direct generalisation of the modal ones, namely as Boolean algebras with a map $\blacksquare$ satisfying the usual modal equation but with values, not in the algebra itself, but in its ideal lattice, that is the set of open elements of its canonical extension (see \cite{Gerhke}). This definition opens the door for a natural definition of validity of a modal formula in a subordination algebra, and therefore for a modal logic investigation.

One last observation. In the dual space of a modal algebra, the accessibility relation is "asymmetrical": the direct image of a clopen subset does not need to be clopen while the inverse one does. Weakening this condition restores the symmetry and if the inverse image may represent a diamond operator, another diamond naturally arises from the direct image. So any subordination algebra (and in particular, any modal algebra) can be endowed with two related modalities. This asks for a comparison of the expressive power of the new obtained languages.

\section{Preliminary dualities}\label{Section1}

In this section, we present and pursue Celani's first results on subordination algebras \cite{Quasimod}. We however adopt the presentation, and notations, of \cite{Sourab} since it does not conveys any specific connotation (neither contact nor modal).

We use classical notations. In particular Boolean operations and constants are denoted by $\wedge, \vee, \neg, 0$ and $1$. For an order $\leq$, we note 
\[ {\uparrow}a  = \lbrace b \mid a \leq b \rbrace \] (and a similar definition for ${\downarrow}a$). For a binary relation $R$ on a set $X$, we note 
\[ R(a) = \lbrace b \in X \mid a \mathrel{R} b \rbrace \text{ and } R^{-1}(a)= \lbrace b \in X \mid b \mathrel{R} a .\rbrace \] And finally, these notations concerning elements of a set are freely extended to subsets of $X$ by 
\[ R[A] = \bigcup \lbrace R(a) \mid a \in A \rbrace \text{ and } R^{-1}[A] = \bigcup \lbrace R^{-1}(a) \mid a \in A \rbrace .  \]

\begin{definition}[\cite{Sourab}]\label{def1} A \textbf{subordination algebra} is a structure $\mathbb{B} = (B,\prec)$ where $B$ is a Boolean algebra and $\prec$ a \textbf{subordination} on $B$, that is a binary relation subject to the following axioms: 
\begin{itemize}
\item[(S1)] $ 0 \prec 0$ and $1 \prec 1$,
\item [(S2)] $ a \prec b,c$ implies $ a \prec b \wedge c$, 
\item [(S3)] $b,c \prec a$ implies $ b \vee c \prec a$,
\item[(S4)] $ a \leq b \prec c \leq d$ implies $a \prec d$.
\end{itemize} 
Equivalently, for each $b \in B$, ${\prec}(b)$ is a filter and ${\prec}^{-1}(b)$ is an ideal of $B$.

An example of subordination algebras is given by the \textbf{contact algebras} (see for instance \cite{Contact}). They are subordination algebras satisfying the axioms (S5) (extensionality), (S6) (reflexivity) and (S7) (symmetry):
\begin{enumerate}[align=left]
\item[(S5)] $ a \neq 0$ implies $b \prec a$ for some $b \neq 0$,
\item[(S6)] $ a \prec b$ implies $a \leq b$,
\item[(S7)] $a \prec b $ implies $\neg b \prec \neg a$.
\end{enumerate}
Moreover, we recall that \textbf{de Vries algebras} (see \cite{Guram1} and \cite{deVries}) are (Boolean) complete contact algebras satisfying axiom (S8) (transitivity): 
\begin{enumerate}[align=left]
\item[(S8)] $ a \prec b$ implies $a \prec c \prec b$ for some $c$.
\end{enumerate}
\end{definition}
Investigations on subordination algebras have to be done in a suitable categorical environment. In the realm of contact algebra, the first morphisms have been introduced in 1962 by de Vries in \cite{deVries} for the particular case of de Vries algebras. From the algebraic point of view, these are very weak morphisms since they are not even Boolean algebra homomorphisms.

Of course, from the model point of view, morphisms should be those Boolean algebra homomorphisms that respect the subordination relation $\prec$. These morphisms have been taken into account in various papers (see for instance \cite{Sourab} or \cite{Contact}), and we shall consider them in this paper too, but under the name of weak morphisms. They will be a common denominator of three other kinds of morphisms we shall introduce now, directly inspired from Celani's q-homomorphisms \cite[Definition 8]{Quasimod}, mainly because they reduce to usual modal algebra homomorphisms when applied to modal algebras, and because their duals are characterised exactly as in the modal case, that is continuous p-morphisms.

Due to the last observation in the introduction, a subordination relation may be considered in two natural ways as a quasi-modal operators (as we shall see in Remark \ref{remark}), leading to four interconnected categories.

\begin{definition}\label{defSubA} Let $\mathbb{B}$, $\mathbb{C}$ be subordination algebras and let $f$ be a map $B \longrightarrow C$. We consider the following axioms: 
\begin{enumerate}[align=left]
\item[(w)] $ a \prec b $ implies $f(a) \prec f(b)$,
\item[($\lozenge$)] $f(a) \prec c$ implies $a \prec b$ and $f(b) \leq c$ for some $b$,
\item[($\blacklozenge$)] $c \prec f(a)$ implies $b \prec a$ and $c \leq f(b)$ for some $b$.
\end{enumerate}
Boolean algebra homomorphisms satisfying (w) will be called \textbf{weak morphisms}, giving rise to the category \textbf{wSub}. The morphisms satisfying (w) and ($\lozenge$) are the \textbf{$\lozenge$-morphisms} or \textbf{white morphisms}, or simply \textbf{morphisms}, and give rise to the white category \textbf{$\lozenge$Sub}, or more simply \textbf{Sub}. Those satisfying (w) and ($\blacklozenge$) are the \textbf{$\blacklozenge$ morphisms} or \textbf{black morphisms} and give rise to the black category \textbf{$\blacklozenge$Sub.} Finally, those satisfying all three axioms are called \textbf{strong morphisms}. They give rise to the strong category \textbf{sSub}.

\end{definition}

Let us now see how these categories  behave toward equivalent presentation of subordination algebras and also toward modal algebras.

\begin{remark}\label{remark}\begin{enumerate}
\item The correspondence between our presentation and Celani's presentation in \cite{Quasimod} goes a follows. A quasi-modal algebra \cite[Definition 1]{Quasimod} is a pair $(B,\Delta)$ where $B$ is a Boolean algebra and $\Delta$ is a quasi-modal operator, that is a map from $B$ into $\mbox{Id}(B)$, the lattice of ideals of $B$, such that:
\begin{enumerate}
\item $\Delta(a\wedge b) = \Delta a \cap \Delta b$,
\item $\Delta 1 = B$.
\end{enumerate}
There exists a bijective correspondence between subordination relations and quasi-modal operators in a Boolean algebra $B$ (cf.~\cite[Theorem 15]{Comp_SUb_Quasi}). This correspondence is established as follows. If $\Delta$ is a quasi-modal operator on $B$ then the relation $\prec_\Delta$ defined by
\[ a \prec_\Delta b \text{ if and only if } b \in \Delta a \]
is a subordination relation. On the other hand, if $\prec$ is a subordination relation, then the operator $\Delta_\prec$ defined by 
\[ \Delta_{\prec}(b) := {\prec^{-1}}(b) \]
is a quasi-modal operator.

Finally, we recall that a q-homomorphism \cite[Definition 8]{Quasimod} is a Boolean homomorphism $f$ such that $\Delta(f(a)) = {\downarrow}f[\Delta a]$. It is then easy to see that a map $f:(B,{\prec}) \rightarrow (C,\prec)$ is a morphism in \textbf{$\blacklozenge$Sub} if and only if it is a q-homomorphism between $(B,\Delta_\prec)$ and $(C,\Delta_\prec)$.
\item The categories \textbf{Sub} and \textbf{$\blacklozenge$Sub} are trivially isomorphic. Indeed, a (white) subordination algebra $(B,\leq, \prec)$ can be associated to a black one: $(B,\leq, \recprec), $ where $a \recprec b$ if and only if $\neg b \prec \neg a$. Hence, any result in  \textbf{Sub} will give an analogue result in \textbf{$\blacklozenge$Sub}.
\item As already mentioned in \cite[p. 384]{Sourab}, modal algebras may be considered as subordination algebras. To set up notations, we recall that a modal algebra is a pair $\mathbb{B}=(B,\lozenge)$ where $B$ is a Boolean algebra and $\lozenge$ is a (diamond) operator on $B$, that is an unary map on $B$ such that $\lozenge (a \vee b) = \lozenge a \vee \lozenge b$  and $\lozenge 0 = 0$. Then, $\mathbb{B}$ is turned in into a subordination algebra by defining 
\begin{enumerate}
\item $a \prec b$ if and only if $\lozenge a \leq b$.
\end{enumerate} 
As such, we obtain a full subcategory \textbf{MA} of \textbf{Sub} (the white category of modal algebras).

Order dually, a modal algebra $(B,\lozenge)$ may be turned into a subordination algebra by defining 
\begin{enumerate}
\setcounter{enumii}{1}
\item $a \prec b$ if and only if $a \leq \blacksquare b$,
\end{enumerate}
(where $\blacksquare = \neg \blacklozenge \neg$). This was the choice of Celani and we obtain a full subcategory \textbf{$\blacklozenge$MA} of \textbf{$\blacklozenge$Sub} (the black category of modal algebras). 

Finally, there is the weak category \textbf{wMA} of modal algebras with Boolean algebra morphisms $f$ satisfying $\lozenge f(a) \leq f(\lozenge a)$, corresponding to axiom (w), whether (a) or (b) is used.
\item In \cite{SantoliThese}, subordination algebras are studied within another language, the language of strict implication, $\rightsquigarrow$, defined on subordination algebras by
\[ a \rightsquigarrow b = \left\lbrace \begin{array}{l}
1 \text{ if } a \prec b \\
0 \text{ if not} 
\end{array} \right. . \]
The advantage is that this language is purely algebraic. However, the class of subordination algebras is not a variety, at least in the usual sense. And, on this precise class, the strict implication language has exactly the same expressive power than the subordination language.  We will see that the subordination language also has a certain - and useful - algebraic flavour: Celani's point of view is that $\prec$ is best being considered as a multivalued operation. At the level of morphisms, it is not difficult to show that the Boolean homomorphisms $h$ respecting $\rightsquigarrow$ are exactly the Boolean morphisms satisfying 
\[ a \prec b \Leftrightarrow h(a) \prec h(b). \]
Moreover, on the class of reflexive subordination algebras, that is subordination algebras satisfying (S6), the $\rightsquigarrow$ morphisms coincide with the weak embeddings as we will define in Definition \ref{def_weak_emb}.
\end{enumerate}
\end{remark}

The topological counterparts of the subordination categories come as no surprise (see \cite{Quasimod} or also \cite{Sourab}) and we have the following.
\begin{definition}
A \textbf{subordination space} is a topological structure $\mathbb{X}=(X,R)$ where $X$ is a Stone space and $R$, the \textbf{accessibility relation}, a closed binary relation on $X$. Such a structure is called descriptive quasi-modal space by Celani in \cite[Definition 10]{Quasimod}. 
\end{definition}

\begin{remark}Throughout the paper, we will freely use the following results, which are folklore: if $\mathbb{X}$ is a subordination space, then $R[C]$ and $R^{-1}[C]$ are closed for every closed subset $C$; if a closed subset $C$ and an open subset $\omega$ satisfy $C \subseteq \omega$, then there exists a clopen subset $O$ such that $C \subseteq O \subseteq \omega$.
\end{remark}
Now, as in Definition \ref{defSubA}, we define four categories whose objects are the subordination spaces.

\begin{definition}
Let $\mathbb{X}, \mathbb{Y}$ be subordination spaces and let $h$ be a map $X \longrightarrow Y$. We consider the following axioms: 
\begin{enumerate}[align=left]
\item[(w)] $x \mathrel{R} y$ implies $h(x) \mathrel{R} h(y)$,
\item[($\lozenge$)]$h(x) \mathrel{R} y$ implies for some $z \in X$, $y = h(z)$ and $x\mathrel{R}z$,
\item[($\blacklozenge$)] $y \mathrel{R} h(x)$ implies for some $z \in X$, $y = h(z)$ and $z \R x$.
\end{enumerate}
As in Definition \ref{defSubA}, continuous maps satisfying (w) are the \textbf{weak morphisms} and give rise to the weak category \textbf{wSubS}. Those satisfying (w) and ($\lozenge$) are the \textbf{white}, or $\lozenge$, \textbf{morphisms}, or simply \textbf{morphisms} and lead to the category \textbf{$\lozenge$SubS}, also denoted \textbf{SubS}. Those satisfying (w) and ($\blacklozenge$) are the \textbf{black morphisms} or \textbf{$\blacklozenge$morphisms} and are the arrows of the category \textbf{$\blacklozenge$SubS}. Finally, those satisfying all three axioms are called \textbf{strong morphisms}. They give rise to the strong category \textbf{sSubS}.

Note that condition (w) together with ($\lozenge$) is widely known in the literature (under the name p-morphism, or also reduction) and that is the reason we choose the white category rather than the black one to be our favourite one.
\end{definition}

The duality between \textbf{Sub} and \textbf{SubS} (and their black and strong counterparts) is then obtained as follows. Let $\mathbb{B}$ be a subordination algebra. Its dual is $\mathbb{X}=(X,R)$ where $X = \mbox{Ult}(B)$ is the ultrafilter space of $B$ with the topology generated by (the clopen sets)
\begin{equation}\label{eq6}
r(b) = \lbrace x \in X \mid x \ni b \rbrace, \ b \in B 
\end{equation} and $R$ is the binary relation on $X$ defined by 
\begin{equation}\label{eq7}		
x\mathrel{R}y \Leftrightarrow {\prec}[y] \subseteq x .
\end{equation}

On the other hand, let $\mathbb{X} = (X,R)$ be a subordination space. Its dual is the subordination space $\mathbb{B} = (B,\prec)$ where $B=\Of(X)$ is the Boolean algebra of all clopen subsets of $X$ and $\prec$ is the binary relation on $B$ defined by 
\begin{equation}\label{eq8}
O \prec U \Leftrightarrow R^{-1}[O] \subseteq U.
\end{equation} 

Now, for the morphisms, if $f : \mathbb{B} \longrightarrow \mathbb{C}$ is a morphism in \textbf{wSub}, then its dual is 
\[ \Ult(f) : \Ult(C) \longrightarrow \Ult(B) : y \longrightarrow f^{-1}(y).\] And, if $h : \mathbb{X} \longrightarrow \mathbb{Y}$ is a morphism in \textbf{wSubS}, then its dual is 
\[ \Of(h) : \Of(Y) \longrightarrow \Of(X) : O \longmapsto h^{-1}(O). \]

With these notations in mind, we have the following theorem.

\begin{theorem}\label{theoremDual}(\cite{Quasimod}) The functors $\Ult$ and $\Of$ establish a dual equivalence between \textbf{wSub} and \textbf{wSubS}. Their restriction also induce a dual equivalence between \textbf{Sub} and \textbf{SubS}, \textbf{$\blacklozenge$Sub} and \textbf{$\blacklozenge$SubS} and finally between \textbf{sSub} and \textbf{sSubS}.
\end{theorem}
\begin{proof}
To extend Celani's proof to the strong category, it is important to observe that the duality between \textbf{$\blacklozenge$Sub} and \textbf{$\blacklozenge$SubS}, established in \cite{Quasimod} thanks to functors slightly different from ours, may also be realised by the functors $\Ult$ and $\Of$, so that the duality \textbf{sSub} - \textbf{sSubS} is just a superposition of the dualities \textbf{Sub} - \textbf{SubS} and \textbf{$\blacklozenge$Sub} - \textbf{$\blacklozenge$SubS}.
\end{proof}

To continue the interpretation of the modal setting to the subordination one, we recall that a modal space is a pair $\mathbb{X}=(X,R)$ where $X$ is a Stone space and $R$ is a closed binary relation on $X$ which is continuous, i.e. for every $O \in \Of(X)$, the set $R^{-1}[O]$ is clopen. We can then consider the full subcategory \textbf{MS} of \textbf{SubS}, whose objects are modal spaces. Let us highlight the fact that when $\mathbb{X}$ is a modal space, then ${\prec}(O)$ is the principal filter generated by $R^{-1}[O]$ (recall \eqref{eq8}). We will soon see that this correspond to a condition to endow subordination algebras with a classical modal structure. In particular, we have that the well-known duality between \textbf{MA} and \textbf{MS} can be  seen as a corollary of the duality between \textbf{Sub} and \textbf{SubS}.

We now give the discrete version of the dualities of Theorem \ref{theoremDual}. The canonical extension of a subordination algebra will then be realised, \`a la Jonsonn-Tarski \cite{JonssonTarski}, as the "composition" of the topological duality with the discrete one.

\begin{definition}\label{defcomplato} Let $\mathbb{B}$ be a subordination algebra. then $\mathbb{B}$ is said to be \textbf{complete atomic} if $B$ is a complete Boolean algebra and $\prec$ is a \textbf{complete subordination}, that is it satisfies, for any family $(b_i \mid i \in I )$ in $B$: 
\begin{enumerate}[align=left]
\item[(S'2)] $ a \prec b_i $ for all $i \in I$ implies $a \prec \bigwedge \lbrace b_i \mid i \in I \rbrace$,
\item[(S'3)] $b_i \prec a$ for all $i \in I$ implies $\bigvee \lbrace b_i \mid i \in I \rbrace \prec a$.
\end{enumerate}

We obtain again four categories whose objects are the complete atomic subordination algebras whose morphisms are complete Boolean homomorphisms. The categories only differ on the choice of the axioms for morphisms as far as the subordination relation is concerned: axiom (w) of Definition \ref{defSubA} will give the subcategory \textbf{wCM} of \textbf{wSub}, axioms (w) and ($\lozenge$) the subcategory \textbf{CM} of \textbf{Sub}, axioms (w) and ($\blacklozenge$) the subcategory \textbf{$\blacklozenge$CM} of \textbf{$\blacklozenge$Sub} and finally, taking all axioms of Definition \ref{defSubA} will give rise to the subcategory \textbf{sCM} of \textbf{sSub}.

Our definition of atomic completeness is given as the natural one in the subordination world, but it fails to reveal an obvious, though important, property of complete atomic subordination algebras: they are in fact - in two different ways - complete modal algebras.

Indeed, (S'2) means that ${\prec}(a)$ is a principal filter so that we may write \[ {\prec}(a) = {\uparrow}(\lozenge a) \] and (S'3) means that the map $a \longmapsto \lozenge a$ is a complete operator, that is, it commutes with supremum of arbitrary families of elements of $B$. Order dually, (S'3) means that ${\prec}^{-1}(a)$ is a principal ideal, so that we may write 
\[ {\prec}^{-1}(a) = {\downarrow}(\blacksquare a) \] and (S'2) means that the map $a \longmapsto \blacksquare a$ is a complete dual operator, that is, it commutes with infimum of arbitrary families of elements of $B$. It follows from this observation that we may consider \textbf{wCM} as a subcategory of \textbf{wMA}, \textbf{CM} as a subcategory of \textbf{MA}, \textbf{$\blacklozenge$CM} as a subcategory of \textbf{$\blacklozenge$MA} and \textbf{sCM} as a subcategory of \textbf{sMA}.
\end{definition}
Since complete atomic subordination algebras and complete atomic modal algebras have been shown to be isomorphic objects, there is no specific discrete duality in the subordination case. We recall the duality for latter uses.

Let $\mathbb{B}=(B,\prec)$ be a complete atomic subordination algebra. Its dual is the Kripke frame $\At(\mathbb{B}) := (X,R)$, where $X = \At(B)$ is the set of atoms of $B$ and the accessibility relation $R$ is given by 
\[ \alpha \mathrel{R} \beta \text{ if } (\forall b \in B)(\beta \prec b \Rightarrow \alpha \leq b). \]

Conversely, if $\mathbb{X}=(X,R)$ is a Kripke frame, its dual is $ \mathcal{P}(\mathbb{X}): =(B,\prec)$ where $B = \mathcal{P}(X)$, the power set of $X$, and for $E,F \subseteq X$, 
\[ E \prec F \text{ if } R^{-1}[E] \subseteq F.
 \]
 
 As might be expected in this paper, we shall consider four categories of Kripke frames: the weak category \textbf{wKF}, whose morphisms satisfy (w), the category \textbf{KF}, whose morphisms satisfies (w) and ($\lozenge$), \textbf{$\blacklozenge$KF}, with (w) and $(\blacklozenge$) and finally the strong category, whose morphisms respect all the axioms.
 
 The object mapping $\mathcal{P}: \mathbb{X} \longmapsto \mathcal{P}(\mathbb{X})$ is extended to a functor by defining $\mathcal{P}(h)=h^{-1}$ and the object mapping $\At : \mathbb{B} \longmapsto \At(\mathbb{B})$ is extended to a functor by defining for $f : \mathbb{B} \longrightarrow \mathbb{C}$ 
 \[ \At(f) : \At(\mathbb{C}) \longrightarrow \At(\mathbb{B}) : \alpha \longmapsto \wedge \lbrace b \in B \mid \alpha \leq f(b) \rbrace. \] 
 
 \begin{theorem} The functors $\At$ and $\mathcal{P}$ establish a dual equivalence between the categories \textbf{wCM} and \textbf{wKF}, \textbf{CM} and \textbf{KF},\textbf{$\blacklozenge$CM} and \textbf{$\blacklozenge$KF} and finally between \textbf{sCM} and \textbf{sKF}.
 \end{theorem}
 \begin{proof}
 This result is folklore \cite{Thomason}.
 \end{proof}
 As we previously stated, we now use the topological and the discrete dualities to obtain a functor from the categories of subordinations algebras to the ones of complete subordination algebras.
 \begin{definition}\label{defext} The composition $\cdot^\delta = \mathcal{P} F \Ult$, where $F$ forget the topology, is called the \textbf{canonical extension functor}, and $\mathbb{B}^\delta = \mathcal{P}(F(\Ult(\mathbb{B})))$ is the \textbf{canonical extension} of $\mathbb{B}$. The functor $\cdot^\delta$ may be considered as a functor from \textbf{wSub} to \textbf{wCM}, from \textbf{Sub} to \textbf{CM}, from \textbf{$\blacklozenge$Sub} to \textbf{$\blacklozenge$CM} or from \textbf{sSub} to \textbf{sCM}.
 
 To complete the picture, we mention that the natural map $r : \mathbb{B} \longrightarrow \mathbb{B}^\delta$, as defined in \eqref{eq6}, is usually not a morphism in \textbf{Sub}, but it is a little more than a weak morphism, as we can see now.
 \end{definition}
 
 \begin{definition}\label{def_weak_emb} Let $\mathbb{B} \text{ and }\mathbb{C}$ be subordination algebras. Then a map $f : \mathbb{B} \longrightarrow \mathbb{C}$ is said to be a \textbf{weak embedding} if it is a one-to-one weak morphism such that $f(b) \prec f(c)$ implies $ b \prec c$.
 \end{definition}
 
 \begin{proposition}\label{Prop1}
 \begin{enumerate}
 \item The natural map $r : \mathbb{B} \longrightarrow \mathbb{B}^\delta$ is a weak embedding.
 \item For each morphism $f$ in \textbf{wSub} from $\mathbb{B}$ into a complete atomic subordination algebra $\mathbb{C}$, there is a unique weak morphism $g: \mathbb{B}^\delta \longrightarrow \mathbb{C}$ such that $g \circ r = f$.
 \end{enumerate}
 \end{proposition}
 \begin{proof}
 The first assertion is clear. Note that the second assertion does not follow from the functioriality of $\cdot^\delta$ because $\mathbb{C}^\delta$ does not necessarily coincide with $\mathbb{C}$ in case the latter is complete atomic. Anyway, the result is well known at the Boolean level. To reach the subordination level, the easiest way is to notice that $g$ is necessarily the dual (in the discrete duality \textbf{wCM} - \textbf{wKF}) of the composition $h = \Ult(f) \circ j$, where $j$ is the natural weak embedding $\At(C) \longrightarrow \Ult(C): \alpha \longmapsto {\uparrow}\alpha $, and $\Ult(f)$ is the dual (in the duality \textbf{wSub} - \textbf{wSubS}) of $f: \mathbb{B} \longrightarrow \mathbb{C}$. So $g$ is weak whenever $f$ is weak since $j$ is weak.
 \end{proof}
 
 \begin{remark}
 Note that the second assertion of \ref{Prop1} does not extend to morphisms in \textbf{Sub}, meaning that $g$ is not necessarily a morphism in \textbf{Sub} in this case, as seen in the case $\mathbb{B} = \mathbb{C}$ and $f$ is the identity.
 \end{remark}
 
 Our definition of canonical extension is not the only possible one. Often, canonical extensions are obtained in a two-step construction. First, consider an already made canonical extension for a reduct of the structure, then device a formula to extend the additional operations to the canonical extension of the reduct. This could have been done here too. In two ways, in fact, subordination algebras are Boolean algebras with an additional relation. The canonical extension $B^\delta$ of a Boolean algebra $B$ is well known: $B^\delta$ is the power set of the dual $X$ of $B$. There are three ways to consider $\prec$.
 
We may consider subordination algebras as Boolean algebras $B$ with the strict implication: though strictly speaking, it is not an operator. Indeed, it is a binary operator as a map $B^\partial \times B \rightarrow B^\partial$ (where $B^\partial$ is the order dual of $B$). But, as shown in \cite{Nonsmooth}, this gives a non smooth extension (see \cite{Nonsmooth} for definitions).

Another way to consider subordination relations is, as it is done in \cite{Quasimod}, as multi-operators $\lozenge : B \longmapsto \mathcal{P}(X) : a \longmapsto {\prec}(a)$.

Finally,  recall that there is a bijective correspondence between the filters of $B$ and the closed element of $B^\delta$ (\cite[Lemma 3.3]{Gerhke}). Therefore, using \cite{Alessandra} denomination, it is possible to consider a subordination relation as the \textbf{slanted} operator 
 \begin{equation} \lozenge : B \longrightarrow B^\delta : a \longmapsto \bigwedge {\prec}(a).\end{equation}
  The formulas giving the sigma extension and the pi extension of slanted operators are described in \cite{Alessandra} and give for an increasing slanted map $f:B \longmapsto B^\delta$: 
 \begin{align*}
 f^\sigma (x) &= \bigvee \lbrace  \bigwedge \lbrace f(b) \mid b \in B, b \geq c \rbrace \mid c \ \text{closed}, c \leq x \rbrace \\
  f^\pi (x) &= \bigwedge  \lbrace \bigvee \lbrace f(b) \mid b \in B, b \leq o \rbrace \mid o \ \text{open}, o \geq x \rbrace
 \end{align*}
 The advantage of the latter point of view is that unary slanted operators are smooth, as shown in the following result.
 
 \begin{proposition}\label{Prop13} The multi-operator associated to the subordination relation in a subordination algebra is smooth. And its extension to the canonical extension of its Boolean part coincides with the canonical extension functor of Definition \ref{defext}.
 \end{proposition}
 \begin{proof}
 The canonical extension functor gives $\mathbb{B}^\delta = (\mathcal{P}(F(\Ult(B)),\prec)$ where $E \prec F$ if $R^{-1}[E] \subseteq F$, so its associated modal operator is $\lozenge^\delta E = R^{-1}[E]$.
 
 We show now that $\lozenge^\sigma E= \lozenge^\pi E= R^{-1}[E]$.
 Indeed, let $\mathcal{O}(X)$ denotes the open subsets of $X$ and $\mathcal{C}(X)$ its closed subsets, we obtain
 \begin{align}
 \lozenge^\pi(E) &= \bigcap \lbrace \bigcup \lbrace R^{-1}[V] \mid \Of(X) \ni V \subseteq O \rbrace \mid  \mathcal{O}(X) \ni O \supseteq E \rbrace \nonumber \\
 &= \bigcap \lbrace R^{-1}[O] \mid \mathcal{O}(X) \ni O \supseteq E \rbrace \nonumber\\
 &\supseteq R^{-1}[E] \nonumber \\ &\supseteq \bigcup \lbrace R^{-1}[F] \mid \mathcal{C}(X) \ni F \subseteq E \rbrace \label{plop2} \\
&= \bigcup \lbrace \bigcap \lbrace R^{-1}[V] \mid \Of(X) \ni V \supseteq F \rbrace \mid \mathcal{C}(X) \ni F  \subseteq E \rbrace = \lozenge^\sigma E \label{plop}\\
 &=  \bigcap \lbrace \bigcup \lbrace R^{-1}[V_{\alpha(F)}] \mid \mathcal{C}(X) \ni F  \subseteq E \rbrace \mid \alpha \in \Lambda \rbrace \nonumber  \\
 & \supseteq \bigcap \lbrace R^{-1}[O] \mid \mathcal{O}(X) \ni O  \supseteq E \rbrace  = \lozenge^\pi(E) \nonumber
 \end{align}
 where $\Lambda = \lbrace \alpha \mid \alpha \text{ choice function } : F \text{ closed } \subseteq E \longmapsto \alpha(F) \text{ clopen } \supseteq F \rbrace$ and the equality \eqref{plop2} $=$ \eqref{plop} is obtained by Esakia's lemma (\cite{Esakia} or \cite{Sahlqvist}).
 \end{proof}
 
 \section{Subordination algebras as models for modal logic}\label{Section2}
 We now define validity of a modal formula in a subordination algebra, taking into account the fact that $\mathbb{B}^\delta$ is a modal algebra. In fact, we noticed that any complete atomic subordination algebra may be endowed with two modalities $\lozenge$ and $\blacklozenge$ (see Definition \ref{defcomplato}). The first as been evaluated in Proposition \ref{Prop13}: if $\mathbb{B}^\delta = \mathcal{P}(\Ult(B))$, then, for every $E \in \mathbb{B}^\delta$
 \begin{equation}\label{eq_def_of_lozenge} \lozenge E = R^{-1}[E].\end{equation}
 The other one is its order dual (see Remark \ref{remark}), and we have, for every $E \in  \mathbb{B}^\delta$,
 \[ \blacklozenge E = R[E]. \]
 
 So, our formulas will be bimodal formulas, that is, terms over the language $(\vee,\wedge,\neg,\top,\bot,\lozenge,\blacklozenge)$. Those formulas not using $\blacklozenge$ will be called simply modal formulas, or white formulas if needed, while formulas not using $\lozenge$ will be called black formulas. Of course, $\mathbb{B}^\delta$, as a bimodal algebra, satisfies formulas not valid in all bimodal algebras since $\lozenge$ and $\blacklozenge$ are both induced by a single accessibility relation. Hence, $\mathbb{B}^\delta$ is what is called in the literature a tense algebra. We recall now basic facts about tense logic and tense algebras.
 
 \begin{definition} A \textbf{tense bimodal logic} (see for instance \cite{Tensedef}) is a provably closed set of bimodal formulas containing all tautologies, the white and black versions of axiom $K$ 
 \begin{enumerate}[align=left]
 \item[($K_\lozenge$)] $\square (\phi \rightarrow \psi) \rightarrow (\square \phi \rightarrow \square \psi)$,
 \item[($K_\blacklozenge$)] $\blacksquare(\phi \rightarrow \psi) \rightarrow (\blacksquare \phi \rightarrow \blacksquare \psi)$,
 \end{enumerate}and the following axioms 
 \begin{enumerate}[align=left]
 \item[$(T_1)$] $\phi \rightarrow \square \blacklozenge \phi$,
 \item[$(T_2)$] $ \blacklozenge \square \phi \rightarrow \phi$.
 \end{enumerate} 
 
 It was shown in \cite{Thomason2} that a bimodal algebra satisfies $T_1$ and $T_2$ (so is a \textbf{tense algebra}) if and only if the accessibility relation associated to $\lozenge$ is the converse of the accessibility relation associated to  $\blacklozenge$. This shows that complete atomic subordination algebras are tense algebras, and their strong category \textbf{sCM} is a subcategory of the category \textbf{TA} of tense algebras (which is also the strong category of modal algebras as initiated in Remark \ref{remark}).
 
 \end{definition}
 
We now arrive at the key definition that will allow us to use subordination algebras as models for modal/tense logic.
 
 \begin{definition}\label{defval} Let $\mathbb{B}$ be a subordination algebra and let $\varphi$ be a bimodal formula. A \textbf{valuation} on $\mathbb{B}$ is a map $v : \Var \longrightarrow \mathbb{B}$. Using the map $r$ of Definition \ref{defext}, the composition $r \circ v$ becomes a map $\Var \longrightarrow \mathbb{B}^\delta$ and as such, extends to a unique homomorphism, also denoted $v$, from the algebra of all bimodal formulas into $\mathbb{B}^\delta$. We say that $\varphi$ is \textbf{valid in $\mathbb{B}$ under the valuation $v$}, denoted \[ \mathbb{B}\models_v \varphi, \] if $v(\varphi) = 1$. Also, as usual, \[ \mathbb{B}\models \varphi \] means $\mathbb{B} \models_v \varphi$ for all valuations $v$. If $\mathcal{K}$ is a class of subordination algebras and $L$ a set of bimodal formulas, $\mathcal{K} \models \varphi$, $\mathbb{B} \models L$ and $\mathcal{K} \models L$ receive their usual meanings. At the dual level, things are rather natural. Indeed, let $\mathbb{X}=(X,R)$ be the dual of $\mathbb{B}$. Then $\mathbb{B} \models \varphi$ is equivalent to $(X,R) \models_v \varphi$  in the classical Kripke semantic meaning, for all valuations $v$ with values in the clopen subsets of $X$.
 \end{definition}
 
 The reason to define validity of modal formulas on subordination algebras via their canonical extension is clear at the dual level. Indeed, in a subordination space $\mathbb{X}$, the valuation of a formula may fail to be a clopen subset of $\mathbb{X}$, that is an element of the associated subordination algebra, as it will be illustrated in Example \ref{exampl}. Hence, we will need the modal machinery available in $\mathbb{B}^\delta$ to properly evaluate formulas.

 \paragraph{}Before turning to completeness results, the following observation is in order. For a set $L$ of formulas, denote by $\mbox{Thm}(L)$ the logic axiomatized by $L$ 
 \[ \mbox{Thm}(L) = \lbrace \varphi \mid L \vdash \varphi\rbrace. \]
 For a class $\mathcal{K}$ of (subordination) algebras, or spaces, denote by $\mbox{Log}(\mathcal{K})$ the logic of $\mathcal{K}$ 
 \[ \mbox{Log}(\mathcal{K}) = \lbrace \varphi \mid \mathcal{K} \models \varphi \rbrace. \] Then, a completeness theorem identifies syntaxic truth with semantic truth, i.e., is of the form $\mbox{Thm}(L) = \mbox{Log}(\mathcal{K})$ for some $L$ and some $\mathcal{K}$. This is an impossible challenge in our case since $\mbox{thm}(L)$ is always a normal modal logic, while this is not always the case for $\text{Log}(\mathcal{K})$. If $\mbox{Log}(\mathcal{K})$ is closed under modus ponens and necessitation, it is not always closed under substitution, as proved by the following example.
 
 \begin{example}\label{exampl} Let $X$ be a Boolean space with an accumulation point $x$. Define $R \subseteq X \times X$ by $y \mathrel{R} z$ if $y=x$ or $y=z$. Then $\mathbb{X}=(X,R)$ is a subordination space whose logic is not a normal modal logic.
 \end{example}
 \begin{proof}
 We first show that $\mathbb{X} \models \varphi \equiv p \rightarrow \lozenge \square p$. Indeed it is not difficult to show that if $O$ is clopen in $X$, then $O = \lozenge \square O$, except when $O \neq \emptyset$ and $O \nni x$, in which case, $\lozenge \square O = O \cup \lbrace x \rbrace$. In all case, $O \subseteq \lozenge \square O$. As a side note, and to come back to the discussion started after Definition \ref{defval}, since $x$ is an accumulation point, it is clear that $O \cup \lbrace x \rbrace$ is not clopen when $x \nin O$. Hence, the valuation of the formula $\lozenge \square p$ is generally not a clopen subset of $X$.
 
 To complete the proof, it suffices to give an instance of $\varphi$ which is not valid in $\mathbb{X}$. Let $\psi \equiv p \wedge \neg \square p$. Then, $\mathbb{X}\not\models \psi \rightarrow \lozenge\square \psi$, as seen when $p$ is evaluated at a proper clopen subset containing $x$ (in this case, $\psi = \lbrace x \rbrace$ and $\lozenge \square \psi = \emptyset$).
 \end{proof}
 
 \begin{nota}\label{notationscheme} We shall see in Theorem \ref{whiteBirk} conditions under which $\mbox{Log}(\mathcal{K})$ is a normal modal logic but an immediate observation is that the substitution rule may be replaced by the use of schemes. To distinguish the formula $\varphi(\overline{p})$ from its \textbf{associated scheme}, we shall write the latter $\varphi(\overline{\psi})$, this expression denotes the collection of formulas $\varphi(\overline{\psi})$ when $\overline{\psi}$ ranges over all modal (or bimodal if needed) tuples of formulas. We arrive at the following completeness results.
  \end{nota}
 \begin{theorem}\label{completeness1} Let $L$ be a set of schemes of modal formulas, and let $\varphi$ be a modal formula. Then the following are equivalent:
 \begin{enumerate}
 \item $L \vdash \varphi$,
 \item for any modal algebra $\mathbb{B}$, $\mathbb{B} \models L$ implies $\mathbb{B} \models \varphi$,
 \item for any subordination algebra $\mathbb{B}$, $\mathbb{B}\models L$ implies $\mathbb{B} \models \varphi$. 
 \end{enumerate}
 \end{theorem}
 \begin{proof}
If $B$ is a modal algebra, then its canonical extension as a modal algebra coincides with the canonical extension of $B$ considered as a subordination algebra as defined in Remark \ref{remark}. So validity in $B$ qua modal algebra is equivalent to validity in $B$ qua subordination algebra. Therefore Item 3 implies Item 2. Moreover Item 1 is equivalent to Item 2 by the standard completeness theorem. Finally, one can prove that Item 1 implies Item 3 by induction on the length of a proof of $\varphi$. Here of course, the fact that $L$ is a set of schemes is essential!
 \end{proof}

\begin{theorem}\label{completeness2} Let $L$ be a set of schemes of bimodal formulas containing the least tense bimodal logic, and let $\varphi$ be a modal formula. Then the following propositions are equivalent:
\begin{enumerate}
\item $L \vdash \varphi$,
\item for any tense algebra $\mathbb{B}$, $\mathbb{B} \models L$ implies $\mathbb{B} \models \varphi$,
\item for any modal algebra $\mathbb{B}$, $\mathbb{B } \models L$ implies $\mathbb{B} \models \varphi$,
\item for any black modal algebra $\mathbb{B}$, $\mathbb{B} \models L$ implies $\mathbb{B} \models \varphi$,
\item for any subordination algebra  $\mathbb{B}$, $\mathbb{B} \models L$ implies $\mathbb{B} \models \varphi$.
\end{enumerate}
\end{theorem}
\begin{proof}
As in Theorem \ref{completeness1}, we have the following chain of implications
\[ \text{Item 2 } \Leftrightarrow \text{ Item 1} \Rightarrow \text{ Item 5 } \Rightarrow \text{ Item 2}. \] 

Let us now prove that Item 3 implies Item 2. If $(B,\lozenge, \blacklozenge)$ is a tense algebra, then its reduct $(B,\lozenge)$ is a modal algebra, and so is a subordination algebra on which bimodal formulas may be evaluated. As in the proof of the previous completeness theorem, validity of a bimodal formula $\varphi$ in $(B,\lozenge, \blacklozenge)$ or in $(B,\lozenge)$ are equivalent. Therefore, Item 3 implies indeed Item 2. Hence, we have 
\[ \text{ Item 3 } \Rightarrow \text{ Item 2 } \Leftrightarrow \text{ Item 5 } \Rightarrow \text{ Item 3}. \]

Finally, $\text{Item 3} \Leftrightarrow \text{ Item 4}$ is obtained by order duality.
\end{proof}

Of course, Theorems \ref{completeness1} and \ref{completeness2} are not really new completeness theorems, and modal algebra remains the favourite algebraic semantic for normal modal logics. But this is a new soundness result, and the theorems might prove to be useful to show non-provability results (in normal logics without the finite model property), since in the search of counterexamples, we may look in the whole class of subordination algebras (much larger, at least in the infinite case, than the ones of modal/tense algebras).

To put it in other words, if there is a subordination algebra $\mathbb{B} \models L$ such that $\mathbb{B} \not\models \varphi$, then there is a modal algebra $\mathbb{C}$ with the same property. We now give a functorial way to move from $\mathbb{B}$ to $\mathbb{C}$, giving an idea of how $\mathbb{C}$ may be more difficult to obtain than $\mathbb{B}$. 

\begin{definition}\label{def_modalisation} Let $\mathbb{B}$ be a subordination algebra and let $r$ be the weak embedding $\mathbb{B} \longrightarrow \mathbb{B}^\delta$. The (white) modal subalgebra of $\mathbb{B}^\delta$ generated by $r(B)$ is called the \textbf{modalisation} of $\mathbb{B}$ and it is denoted by $\mathbb{B}^m$ and we have the following result.
\end{definition}

\begin{proposition} The object mapping $\mathbb{B} \longrightarrow \mathbb{B}^m$ can be extended to a covariant functor $\cdot^m : \textbf{Sub} \longrightarrow \textbf{MA}$. The natural map $r : \mathbb{B} \longrightarrow \mathbb{B}^m$ is a weak embedding.
\end{proposition}
\begin{proof}
Suppose $f : \mathbb{B} \longrightarrow \mathbb{C}$ is a morphism in \textbf{Sub}. By the canonical extension functor, $f$ lifts to $f^\delta : \mathbb{B}^\delta \longrightarrow \mathbb{C}^\delta$ in \textbf{CM}. Let $f^m$ be the restriction of $f^\delta$ to $\mathbb{B}^m$. It suffices now to show that $f^m$ takes value into $\mathbb{C}^m$.

If $b \in \mathbb{B}^m$, there are $b_1,...,b_n \in B$ and  a modal formula $\varphi$ with $b = \varphi(b_1,...,b_m)$. Then, 
\[ f^\delta(b) = f^\delta(\varphi(b_1,...,b_n)) = \varphi(f^\delta(b_1),...,f^\delta(b_n)) \in  \mathbb{C}^m\] as required.
\end{proof}

Of course, next to the \textbf{modalisation functor}, there is the \textbf{black modalisation functor} $\cdot^{blm}: \textbf{$\blacklozenge$Sub} \longrightarrow \textbf{$\blacklozenge$MA}$ ($\mathbb{B}^{blm}$ is the black modal algebra generated by $r(B)$) and the \textbf{bimodalisation functor} $\cdot^{bim} : \textbf{sSub} \longrightarrow \textbf{TA}$ ($\mathbb{B}^{bim}$ is the least tense algebra generated by $r(B)$).

The relevance of these concepts is that they preserve, not validity of formulas, but validity of schemes, as shown in the following result. We give the results in case of modalisation, but of course, there is an analogue result for black and bi-modalisation.
\begin{proposition}\label{modal1} Let $\mathbb{B}$ be a subordination algebra. Then for any scheme $\varphi(\overline{\psi})$ of modal formulas, 
\[ \mathbb{B } \models \varphi(\overline{\psi}) \Leftrightarrow \mathbb{B}^m \models \varphi(\overline{\psi}). \]
\end{proposition}
\begin{proof}
Since $\mathbb{B}^m$ is submodal algebra of $\mathbb{B}^\delta$, the if part directly follows from the definition of validity in \ref{defval}.

Suppose $\overline{\psi} =(\psi_1,...,\psi_n)$ and $\mathbb{B} \models \varphi(\overline{\psi})$. We have to show $\mathbb{B}^m \models \varphi(\overline{\psi})$. Let $v$ be a valuation $\Var \longrightarrow B^m$. There are $b_1,...,b_r \in B$ and formulas $\varphi_1,...,\varphi_n$ in the variables $q_1,...,q_r$ such that $v(\psi_i) = \varphi_i(b_1,...,b_r)$ for all $i$. Let $v'$ be the valuation $\Var \longrightarrow B$ such that $v'(q_i) = b_i$ for all $i$. Then, $\mathbb{B} \models_{v'} \varphi(\varphi_1,...,\varphi_n)$ and it follows 
\begin{align*}
1&= v'(\varphi(\varphi_1,...,\varphi_n)) \\
&= \varphi(v'(\varphi_1),....,v'(\varphi_n)\\
&= \varphi(\varphi_1(b_1,...,b_r),...,\varphi_n(b_1,...,b_r)) \\
&= \varphi(v(\psi_1),...,v(\psi_n))\\
&= v(\varphi(\overline{\psi})),
\end{align*} as required.
\end{proof}

Note that Example \ref{exampl} shows that a single modal formula may fail to be preserved by modalisation: if $\mathbb{B}$ is the dual of $\mathbb{X}$ as described in \ref{exampl}, then $\mathbb{B} \models p \rightarrow \lozenge \square p$ while $\mathbb{B}^m \not \models p \rightarrow \lozenge \square p$. Otherwise, by modal logic, $\mathbb{B}^m \models \psi \rightarrow \lozenge \square \psi$ for all modal formula $\psi$ and thus, by Proposition \ref{modal1}, $\mathbb{B} \models \psi \rightarrow \lozenge \square \psi$, which is not the case.

We now continue our discussion begun just before Definition \ref{def_modalisation}. The distance between $\mathbb{B}$ and $\mathbb{B}^m$, which is also the measure of how the accessibility relation $R$ fail to be continuous, is fully realised at the dual level. 

\begin{proposition}\label{prop_dual_modalisation} Let $\mathbb{B}=(B,\prec)$ be a subordination algebra whose dual space is $\mathbb{X}=(X,R)$.
\begin{enumerate}
\item The topological dual of $\mathbb{B}^\delta$ is given by $\mathbb{X}^\delta=(\beta(X), \overline{R})$, where $\beta(X)$ is the Stone-\v{C}ech compactification of $X$ endowed with the discrete topology and $\overline{R}$ is the closure of $R$ in $\beta(X)$.
\item The dual of $\mathbb{B}^m$ is the quotient of $\mathbb{X}^\delta$ by the equivalence relation $\sim$ defined by $u \sim v$ if and only if $\varphi(\underline{b}) \in u $ if and only if $\varphi(\underline{b}) \in v$ for any modal formula $\varphi$ and any tuple $\underline{b} \in B^m$.
\end{enumerate}
\end{proposition}
\begin{proof}
Let us denote by $\iota$ the embedding 
\begin{equation}\label{eq_def_iota} \iota : X \hookrightarrow \Ult(\mathcal{P}(X)) \cong \beta(X) : x \mapsto \lbrace E \in \mathcal{P}(X) : x \in E \rbrace, \end{equation}
and by $R^\beta$ the accessibility relation associated to $\mathbb{B}^\delta$. We then have the following diagram. 

\begin{center}
\begin{tikzpicture}
\matrix (m) [matrix of math nodes,row sep=3em, column sep=3em, minimum width=2em]
{(B,\prec) & (B^\delta,\lozenge)  \\ (X,R)  & (\Ult(\mathcal{P}(X)),R^\beta) \\};
\path[-stealth]
(m-1-1)  edge  (m-2-1) 
(m-2-1) edge   (m-1-1)
(m-1-1) edge node[above]{$r$} (m-1-2)
(m-2-1) edge node [below] {$\iota$} (m-2-2)
(m-1-2) edge (m-2-2)
(m-2-2) edge (m-1-2);
\end{tikzpicture}
\end{center}
\begin{enumerate} 
\item The result for the underlying structure of $\mathbb{X}^\delta$ is folklore. For the relation, we have first that $\overline{R} \subseteq R^\beta$. 

Indeed, for $x,y \in X$, we have 
\begin{center}
\begin{tabular}{rlll}
& $\iota(x) \mathrel{R}^\beta \iota(y)$ & \\
$\Leftrightarrow$ & $\lozenge[\iota(y)] \subseteq \iota(x)$ && By definition of standard modal dual \\
$\Leftrightarrow$ & $ F \ni y  \Rightarrow \lozenge F \ni x$ && By \eqref{eq_def_iota} \\ $\Leftrightarrow$ & $F \ni y \Rightarrow R^{-1}[F] \ni x$ &&By definition of ${\lozenge}$ (see \eqref{eq_def_of_lozenge}).
\end{tabular}
\end{center}
Therefore, $x \mathrel{R} y$ clearly implies $\iota(x) \mathrel{R}^\beta \iota(y)$ and, since $R^\beta$ is closed, it follows that $\overline{R} \subseteq R^\beta$.

On the other hand, let $u,v \in \Ult(\mathcal{P}(X))$ be such that $(u,v) \nin \overline{R}$. Then, since $\overline{R}$ is closed, there exist $E_1, E_2 \in \mathcal{P}(X)$ such that $E_1 \in u$, $E_2 \in v$ and  $(E_1 \times E_2) \cap R = \emptyset$, which is equivalent to 
\[ \lozenge E_2 = R^{-1}[E_2] \subseteq E_1^c. \] 
Hence, we have $E_2 \in v$ and $\lozenge E_2 \nin u$ (indeed $u$ is an ultrafilter containing $E_1$). It follows that $(u,v) \nin R^\beta$ and the proof is completed. 
\item Let us denote by $\mathbb{X}^m$ the Stone dual of $\mathbb{B}^m$. Since $\id : \mathbb{B}^m \rightarrow \mathbb{B}^\delta$ is blatantly a one-to-one morphism, we know that $\Ult(\id) : \mathbb{X}^\delta \rightarrow \mathbb{X}^m$ is an onto function. Since, for $u,v \in \mathbb{X}^\delta$, we have
\[ \Ult(\id)(u) = \Ult(\id)(v) \text{ if and only if } \mathbb{B}^m \cap u = \mathbb{B}^m \cap v, \] the conclusion follows from the definition of $\mathbb{B}^m$.
\end{enumerate}
\end{proof}
\section{Universal algebraic approach of subordination algebras}\label{Section_universal algebra}

In this section, we look for constructions that preserve the validity of formulas. Since the validity of a formula $\varphi$ is equivalent to the validity of the equation $\varphi = 1$, we naturally turn to the HSP-theorem of Birkhoff and look for a universal algebraic treatment of subordination algebras, in the spirit of \cite{IntroUniv}. This work has been undertaken by Celani in \cite{qcong} and by Celani and Castro in \cite{Castro} as far as subalgebras and quotients are concerned. We present and continue their work using a slightly different terminology.

\begin{definition}\label{defcong} Let $\mathbb{B}$ be a subordination algebra. We say that a binary relation $\theta$ on $B$ is a \textbf{white} (resp. \textbf{black}, \textbf{strong}) \textbf{congruence} if it is the kernel of a white (resp. black, strong) morphism, that is there is a morphism $f: \mathbb{B} \longrightarrow  \mathbb{C}$ such that
\[ \theta = \ker(f) = \lbrace (a,b) \in B^2 \mid f(a) = f(b) \rbrace. \]
Of course, a white congruence is necessarily a (Boolean) congruence. We recall that congruences on a Boolean algebra $B$ are in bijective correspondence with ideals and filters of $B$ respectively via the assignments 
\[ \theta \longmapsto I_\theta := \lbrace a \in B \mid a \mathrel{\theta} 0 \rbrace \] and
\[ \theta \longmapsto F_\theta := \lbrace a \in B \mid a \mathrel{\theta} 1 \rbrace. \]
Also, we recall that $\theta, I_\theta $ and $F_\theta$ are linked by the following  formulas: $a \mathrel{\theta} b$ if and only if there is $i \in I_\theta$ with $a\vee i = b \vee i$, if and only if there is $f\in F_\theta$ with $a \wedge f = b \wedge f$. Finally, we recall that an ideal $I$ is \textbf{round} (cf.~\cite{Guram1}) if $a \in I$ implies $a \prec b$ for some $b\in I$, and that a filter $F$ is \textbf{round} if $a \in F$ implies $b \prec a$ for some $b \in F$. 
\end{definition}

Let us compare with the congruences of Celani. If $(B,\Delta)$ is a quasi-modal algebra, a quasi-modal congruence \cite[Definition 10]{qcong} is a Boolean congruence $\theta$ such that $a \mathrel{\theta} b$ and $c \in \Delta a$ implies $c \mathrel{\theta} d$ for some $d \in \Delta b$. Implicitly in \cite{qcong}, we have quasi-modal congruences correspond to the kernels of quasi-modal morphisms. Hence, by Remark \ref{remark}, we have that a Boolean congruence $\theta$ is a quasi-modal congruence on $(B,\Delta)$ if and only if it is a black congruence of $(B,\prec_\Delta)$.

In fact, it is a matter of convenience to define congruences by some substitution property or as a kernel of morphisms (we choose the latter because it applies in any algebraic category). We give here the adaptation of \cite[Theorems 16 and 17]{qcong} to our white category.

\begin{proposition}\label{prop3.7} Let $\mathbb{B}=(B,\prec)$ be a subordination algebra and let $\theta$ be a congruence on $B$. With the notations of Definition \ref{defcong}, the following are equivalent: 
\begin{enumerate}
\item $\theta$ is a white congruence,
\item $\theta$ satisfies $a \mathrel{\theta}b \prec c $ implies $ a \prec d \mathrel{\theta} c $ for some $d$,
\item $I_\theta$ is a round ideal,
\item $F_\theta$ satisfies: $a \in F_\theta$ implies $\neg a \prec \neg b$ for some $b \in F_\theta$.
\end{enumerate}
\end{proposition}

As a direct consequence of Proposition \ref{prop3.7}, we have that white congruences are in bijective correspondence with round ideals and what could be call negative round filters. Remark that for subordination algebras satisfying S7, negative round filters and round filters are equivalent notions.

The dual objects of white congruences are easy to characterise: they are the closed subsets $C$ of the dual space which are \textbf{$R$-increasing}, that is, such that $x \in C$ and $x  \mathrel{R} y$ imply $y \in C$. The correspondence, already described in \cite{qcong}, is obviously the restriction to white congruences of the correspondence between (Boolean) congruences (or ideals) and closed subsets of the dual: 
\[ I \text{ ideal of } B \longmapsto C = \lbrace x \in \Ult(B) \mid x \cap I = \emptyset \rbrace. \] As a corollary, we have the following result which is a small improvement of \cite[Lemma 9]{qcong}.

\begin{corollary}\label{prop3.9} Let $\mathbb{B}$ be a subordination algebra and let us denote by $\Con(\mathbb{B})$ the ordered set of all white congruences on $\mathbb{B}$. Then $\Con(\mathbb{B})$  is a frame (complete Heyting algebra) in which finite meet is intersection while arbitrary joint is joint in the equivalence lattice.
\end{corollary}
\begin{proof}
By Proposition \ref{prop3.7}, it suffices to argue on round ideals. Let $I$ and $J$ be round ideals. If $ a \in I \cap J$, there are $b \in I$ and $c \in J$ such that $a \prec b$ and $a \prec c$. Then $ a \prec b \wedge c$ and $b\wedge c \in I \cap J$ (of course this argument does not work for infinitely many round ideals).

Let now $(I_l \mid l \in \Lambda)$  be round ideals and  $I$ its joint in the ideal lattice. If $a \in I$, there are $l_1,\ldots,l_n \in \Lambda$ and $a_1 \in I_{l_1},\ldots,a_n \in I_{l_n}$ with $a \leq a_1 \vee \cdots \vee a_n$. And there are $b_1 \in I_{l_1},\ldots,b_n \in I_{l_n}$ with $a_1 \prec b_1, \ldots,a_n \prec b_n$. It follows that $ a \prec b = b_1 \vee \cdots \vee b_n \in I$.

It follows that $\Con(\mathbb{B})$ is a frame, being a subset of the ideal lattice of $B$, closed under finite meets and arbitrary joints.
\end{proof}

To continue to show that subordination algebras are more on the side of universal algebra than relational structures, it is time to adapt the three classical isomorphism theorems. We omit the classical proofs.

\begin{proposition} Suppose $\mathbb{B}=(B,\prec)$ is a subordination algebra and $\theta$ is a white congruence on $\mathbb{B}$. The structure $\mathbb{B/_\theta} =(B/_\theta, \prec^\theta)$ with
\begin{equation}\label{eq_def_prec_theta} a^\theta \prec^\theta b^\theta \text{ if } (\exists c )(  a \prec c \mathrel{\theta}b) \end{equation} is a subordination algebra such that the canonical projection $\pi : B \longrightarrow B/_\theta$ is a (white) morphism.
\end{proposition}
\begin{proof}
We prove that the relation $\prec^\theta$ is well defined and leave the remaining verifications (namely, $\prec^\theta$ is a subordination and $\pi$ is a morphism) to the reader. Suppose that $a_1 \mathrel{\theta} a_2$, $b_1 \mathrel{\theta} b_2$ and $a_1^\theta \prec^\theta b_1^\theta$. It follows that 
\[ a_2 \mathrel{\theta} a_1 \prec c \mathrel{\theta} b_1 \mathrel{\theta} b_2 \] for some $c \in B$. Now, since $\theta$ is a white congruence, Proposition \ref{prop3.7} implies that we have $ a_2 \prec d \mathrel{\theta} b_2$ for some $d \in B$. Hence, we have $a_2^\theta \prec^\theta b_2^\theta$, as required.
\end{proof}

\begin{proposition}\label{prop3.10}Suppose $\mathbb{A}$ and $\mathbb{B}$ are subordination algebras.
\begin{enumerate}
\item Let $f : \mathbb{A} \longrightarrow \mathbb{B}$ be a morphism and let $\theta$ be a white congruence on $\mathbb{A}$. Then, there is a morphism $g : \mathbb{A/\theta} \longrightarrow \mathbb{B}$ such that $ g \circ \pi = f$ if and only if $\theta \subseteq \ker(f)$. In particular, $\mathbb{A}/\ker(f)$ is isomorphic with the range of $f$, which appears to be a subordination subalgebra of $\mathbb{B}$ (see \cite[Theorem 18]{qcong}).
\item Suppose $\mathbb{A}$ is a subalgebra of $\mathbb{B}$ and let $\theta$ be a white congruence on $\mathbb{B}$. Then the restriction of $\theta$ to $\mathbb{A}$ is a white congruence on $\mathbb{A}$, the saturation 
\[ \mathbb{A}^\theta := \lbrace b \in \mathbb{B} \mid \exists \ a\in \mathbb{A}, \ a\mathrel{\theta} b \rbrace \] is a subalgebra of $\mathbb{B}$ and $\mathbb{A/\theta|_A}$ is isomorphic with $\mathbb{A^\theta/\theta}$.
\item Let $\theta$ be a white congruence on $\mathbb{B}$. Then the white congruences lattice of $\mathbb{B/\theta}$ is isomorphic with the principal filter of $\Con(\mathbb{B})$ generated by $\theta$.
\end{enumerate}
\end{proposition}

Of course, analogue to Propositions \ref{prop3.7} to \ref{prop3.10}, which are relative to \textbf{Sub}, there are corresponding results in the categories \textbf{$\blacklozenge$Sub} and \textbf{sSub}. For the sake of completeness, let us note that black congruences are the Boolean congruences which are characterised by round filters as a 1-kernels. Moreover, strong congruences are Boolean congruences whose 0-kernels are round ideals and 1-kernels are round filters.

\begin{proposition}\label{quotient} If $\varphi$ is a modal formula (resp. a black modal formula or a bimodal formula), $\mathbb{B}$ a subordination algebra and $\theta$ a white congruence (resp. a black congruence or a strong congruence), then 
\[ \mathbb{B} \models \varphi \ \Rightarrow \ \mathbb{B/_\theta} \models \varphi. \] In other words, validity of a formula is preserved by morphic image of its language.
\end{proposition}
\begin{proof}
It suffices to examine the white language.
Let $f$ be an onto morphism (in \textbf{Sub}) from $\mathbb{B}$ to $\mathbb{C}$. Suppose $\mathbb{B} \models \varphi$. We want to prove that $\mathbb{C} \models \varphi$. Choose a valuation $v : \Var \longrightarrow C$  and let $v_1$ be any valuation $\Var \longrightarrow B$ such that $f \circ v_1 = v$. Then $v_1$ extend to a modal homomorphism $v_1 : \mbox{Form} \longrightarrow \mathbb{B}^\delta$ such that $v_1(\varphi) = 1$. Now, applying the canonical extension functor, we have a morphism $f^\delta : \mathbb{B}^\delta \longrightarrow \mathbb{C}^\delta$, so $f^\delta \circ v_1$ is the extension to $\mbox{Form} \longrightarrow \mathbb{C}^\delta$ of $v$. And we have $v(\varphi) = f^\delta(v_1(\varphi)) = 1$ as required. 
\end{proof}

We now turn to subalgebras (subobjects in \textbf{Sub}), a topic examined in \cite{Castro} for lattices in the quasi-modal setting. Our condition \eqref{eq_eq_eq_8}, or more precisely, its black version, corresponds to the quasi-modal condition required by Castro and Celani in \cite[Definition 21]{Castro}.

\begin{definition}\label{def_subsubalg} Let $\mathbb{B}$ be a subordination algebra and $A \subseteq B$. Then $A$ is a \textbf{subalgebra} of $\mathbb{B}$ if $A$ is a (Boolean) subalgebra of $B$ and 
\begin{equation}\label{eq_eq_eq_8} \text{ if } a \in A, \ b \in B \text{ and } a \prec b, \ \text{there is } a \in A \text{ with } a \prec c \leq b.
\end{equation}
Of course the associate subordination algebra is $\mathbb{A}=(A, \prec)$ where $\prec$ is the restriction to $A$ of the subordination relation of $\mathbb{B}$. We also say that $\mathbb{A}$ is a subalgebra of $\mathbb{B}$. This is the appropriate concept since $\mathbb{A}$ is a subalgebra if and only if the inclusion mapping $i:\mathbb{A} \longrightarrow \mathbb{B}$ is a morphism. And moreover, if $\mathbb{B}$ is a modal algebra, the concepts of subalgebra of $\mathbb{B}$, qua modal algebra qua subordination algebra, coincide.
\end{definition}

However, this concept of subalgebra lacks many usual properties of subalgebras of universal algebra. We shall see in \ref{plop1} that the intersection of two subordination subalgebras may not be a subordination subalgebra. So it is even not clear that the set $\mbox{Sub}(\mathbb{B})$ of all subalgebras of $\mathbb{B}$ forms a lattice. The only obvious positive result we have in this direction is the following.
 
 \begin{proposition} If $( A_i \mid i \in I )$ is a directed family of subalgebras of $\mathbb{B}$, so is its union. Hence, $\Sub (\mathbb{B})$ is a dcpo in the sense of \cite{Compendium}. 
 \end{proposition}
 
 We now characterise the notion dual to subalgebras. One again, we have a restriction of a known correspondence in Stone duality: the correspondence between Boolean subalgebras and congruences of Stone spaces. 
 
 \begin{definition}\label{defcongX}Let $\mathbb{X} = (X,R)$ be a subordination space. A \textbf{congruence on $X$} (or more precisely a \textbf{white congruence}) is an equivalence relation $\Theta$ on $X$ such that:
 \begin{enumerate}
 \item $\Theta$ is a congruence of Stone spaces, that is non equivalent points can be separated by clopen $\Theta$-saturated sets (we recall that a subset $S$ is $\Theta$-saturated if $x \in S$ and $x \mathrel{\Theta}y$ implies $y \in S$, that is $S$ is a union of equivalence classes),
 \item  for $x,y,z \in X$ \begin{equation}\label{neweq} x \mathrel{\Theta}y \mathrel{R} z \Rightarrow \exists \ u \in X : x \mathrel{R}u \mathrel{\Theta}z.\end{equation} \end{enumerate}
 
 The condition \eqref{neweq}, or once again, more precisely, the black version given in \eqref{eq_black_cong}, corresponds to the $\Delta$-part of condition 2) in  \cite[Theorem 24]{Castro}. 
 
 \end{definition}
 
Here again this is the appropriate concept since an equivalence $\Theta$ on $\mathbb{X}$ is a white congruence if and only if there is a, necessarily unique, subordination structure on $X/\Theta$ such that the natural projection is a morphism. For $x,y \in X$, we have $x^\Theta \mathrel{ R^\Theta } y^\Theta$ if there is $y' \in X$ such that $x\mathrel{R}y'\mathrel{\Theta}y$, or equivalently, if for any $x' \mathrel{\Theta}x$, there is $y' \in X$ such that $x' \mathrel{R} y' \mathrel{\Theta}y$.

Let us denote by $\Con(\mathbb{X})$ the ordered set of white congruences on $\mathbb{X}$. The duality $\textbf{Sub} \leftrightarrow \textbf{SubS}$ of Theorem \ref{theoremDual} exchanges one-to-one maps with onto ones, so exchanges subalgebras with congruences, and we have the following result.

 \begin{proposition}[\cite{Castro}]\label{prop3.5} Let $\mathbb{B}$ be a subordination algebra and let $\mathbb{X}$ be its dual. Then $\Sub(\mathbb{B})$ is anti-isomorphic (that is isomorphic to the order dual of) to $\Con(\mathbb{X})$.
 \end{proposition}
 \begin{proof}
 The proof is a direct byproduct of the duality. Recall the well-known anti-isomorphism between Boolean subalgebras $A$ of $B$ and Boolean space congruences (see for instance \cite[Chapter 37]{IntroBool}) given by
 \begin{enumerate}
 \item $x  \mathrel{\Theta} y$ if and only if $x \cap A = y \cap A$, and 
 \item $A = \lbrace O \in \Of(X) \mid O \text{ is $\Theta$-saturated } \rbrace$
 \end{enumerate}
 Then, $A$ is a (subordination) algebra of $B$ if and only if $\Theta$ is a white congruence, as proved in \cite[Theorem 24]{Castro} for the black setting.
 \end{proof}

 \begin{example}\label{plop1} Let $\omega$ be the set of natural numbers and $\omega^+$ be its successor ordinal (topologically, $\omega^+$ is the Alexandroff compactification of $\omega$). We consider $\omega^+$ as a subordination space with the relation $R$ defined by $x  \mathrel{R} y$ if $y = x$ or $y = \omega$ or $x =\omega$. It is not difficult to show that an equivalence is a congruence if and only if all its classes are closed and the class containing $\omega$ is either $\lbrace \omega \rbrace$ or $\omega^+$. Let $\Theta$ be the equivalence whose classes are $\lbrace 2i,2i+1 \rbrace$ for $i \in \omega$ and $\lbrace \omega \rbrace$; let $\xi$ be the equivalence whose classes are $\lbrace 0,1 \rbrace$, $\lbrace 2 \rbrace$,$\lbrace 2i+3,2i+4\rbrace$ for $i \in \omega$ and $\lbrace \omega \rbrace$. Then, the supremum of $\Theta$ and $\xi$ in the lattice of Boolean congruences of $\omega^+$ has two classes $\lbrace 0,1 \rbrace$ and its complement, hence it is not a congruence. By Proposition \ref{prop3.5}, this shows that the intersection of two subalgebras of a subordination algebra may fail to be a subalgebra.
 \end{example}
 
 As was the case for congruences, there are \textbf{black subalgebras} $A$ of $\mathbb{B}$ ($ a \prec b$ and $a \in A$ implies $b \leq c \prec a $ for some $c \in A$) with dual black congruences on $\mathbb{X}$ satisfying the first condition of Definition \ref{defcongX} and 
 \begin{enumerate}
 \item[2'.] For $x,y,z \in X$
 \begin{equation}\label{eq_black_cong} x \mathrel{R} y \mathrel{\Theta} z \Rightarrow \exists \ u \in X : x \mathrel{\Theta} u \mathrel{R} z ; \end{equation}
  \end{enumerate} and there are \textbf{strong subalgebras} (both black and white subalgebras) and strong congruences of subordination spaces (both white and black congruences).

  \begin{proposition}\label{subalgebra} If $\varphi$ is a modal formula (resp. a black modal formula or a bimodal formula), $\mathbb{B}$ a subordination algebra and $\mathbb{A}$ a subalgebra (resp. a black subalgebra or a strong subalgebra) then 
  \[ \mathbb{B} \models \varphi \Rightarrow \mathbb{A} \models \varphi. \]
  \end{proposition}
  \begin{proof}
  As in Proposition \ref{quotient}, it suffices to examine the white language. The inclusion morphism $i : \mathbb{A} \longrightarrow \mathbb{B}$ lifts to a \textbf{CM}-morphism $i^\delta : \mathbb{A}^\delta \longrightarrow \mathbb{B}^\delta$ which is one-to-one since $\cdot^\delta = \mathcal{P}F\Ult$ and both $\Ult$ and $\mathcal{P}$ exchange onto with one-to-one.
  
  We suppose $\mathbb{B} \models \varphi$ and want to prove $\mathbb{A} \models \varphi$. So let $v$ be a valuation $\Var \longrightarrow \mathbb{A}$. It extends to a modal homomorphism $v: \mbox{Form} \longrightarrow \mathbb{A}^\delta$ and therefore the map $i^\delta v : \mbox{Form} \longrightarrow \mathbb{B}^\delta$ is the homomorphic extension of the valuation $v$ considered as a valuation on $\mathbb{B}$. Since $\mathbb{B} \models \varphi$, we have $i^\delta(v(\varphi)) = 1$, whence $v(\varphi) = 1$, as required.
  \end{proof}

  We now turn to products, a topic which was not examined previously. We concentrate on Cartesian products (defined pointwise).
  
  \begin{proposition}\label{product} Let $(\mathbb{A}_j \mid j \in J)$ be a family of subordination algebras and let $\mathbb{P} = \prod_{j \in J}\mathbb{A}_j$ be its Cartesian product. Then 
  \begin{enumerate}
  \item $\mathbb{P}$ is the categorical product in the weak category \textbf{wSub},
  \item the projections are morphisms in the strong category \textbf{sSub},
  \item a finite product is a categorical one in the categories \textbf{Sub}, \textbf{$\blacklozenge$Sub} and \textbf{sSub}.
  \end{enumerate}
  \end{proposition}
  \begin{proof}
  Assertions $1.$ and $2.$ follow from direct calculations and we prove assertion $3.$ for the category \textbf{Sub}.
  
  Suppose $\mathbb{B}=(B,\prec)$ is a subordination algebra and $f_j$ are morphisms $\mathbb{B} \longrightarrow \mathbb{A}_j$ for $j \in J$. The Cartesian product $f : b \longmapsto (f_j(b) \mid j \in J)$ is a weak morphism and we show it is a morphism in \textbf{Sub}, that is, we prove axiom $(\lozenge)$ of \ref{defSubA}. So let $f(b) \prec c$, we have $f_j(b) \prec c_j$ for all $j$ and there are $d_j$ in $B$ with $b \prec d_j$ and $f_j(d_j) \leq c_j$. Since, $J$ is finite, $d = \wedge_{j \in J}d_j$ exists in $B$ and we have $b \prec d$ and $f(d) \leq c$, as required.
  \end{proof}
  
  The following result shows that for a Cartesian product, being a categorical product - in \textbf{Sub}, \textbf{$\blacklozenge$Sub} or \textbf{sSub} - is a very strong property, except for finite products and the trivial case where the $\mathbb{A}_j$'s are modal algebras.
  
  \begin{proposition} If $\mathbb{A}=(A, \prec)$ is a subordination algebra which is not a modal algebra, some Cartesian power of $\mathbb{A}$ is not a power in the category \textbf{Sub}.
\end{proposition}
\begin{proof}
Let $a \in A$ be such that ${\prec}(a)$ is not a principal filter, and let $(b_j \mid j \in J)$ be another notation for the set ${\prec}(a)$. Then $A^J$ is not a product in \textbf{Sub}. Too see this, let $A_j = A$ for any $j \in J$ (so that $A^J = \prod_{j \in J} A_j$) and let $f_j : A \longrightarrow A_j$ be the identity map for any $j \in J$. The Cartesian product $f$ of the $f_j$ is the diagonal map 
\[ a \longmapsto (a_j \mid j \in J) \] with $a_j = a$ for all $j \in J$. This is not a morphism: let $c \in A^J$ be defined by $c_j = b_j$ for $j \in J$. Then $f(a) \prec c$ but there is no $b \in A$ such that we have $a \prec b$ and $f(b) \leq c$. 
\end{proof}

The following results examine validity of formulas in products.

\begin{proposition} Let $\varphi$ be a bimodal formula. If $(\mathbb{A}_j \mid j \in J)$ is a finite family of subordination algebras then 
\[ \prod_{j \in J} \mathbb{A}_j \models \varphi \Leftrightarrow \mathbb{A}_j \models \varphi \ \forall j. \]
\end{proposition}
\begin{proof}
The if part follows from \ref{quotient} and \ref{product}.

Suppose now $\mathbb{A}_j \models \varphi$ for all $j$. We have to prove $\prod \mathbb{A}_j \models \varphi$. So let $v$ be a valuation on $\prod \mathbb{A}_j$. Then, if $p_j$ denotes the projection from $\prod A_j$ into $A_j$, $p_j \circ v$ is a valuation on $\mathbb{A}_j$ and $p_j(v(\varphi)) = 1$ in $\mathbb{A}_j^\delta$. It follows $v(\varphi) = 1$ in $\prod \mathbb{A}_j^\delta$ and since 
\[ \prod \mathbb{A}_j^\delta \cong (\prod \mathbb{A}_j)^\delta, \] $v(\varphi) = 1$ in $(\prod \mathbb{A}_j)^\delta$, that is $\prod \mathbb{A}_j \models \varphi$.
\end{proof}

For infinite products, the proof does not work since $(\prod \mathbb{A}_j)^\delta$ is not isomorphic with $\prod \mathbb{A}_j^\delta$, in general.

\begin{lemma}\label{prop2.25} Let $\mathbb{A}_j$ $(j \in J)$ be a subordination algebra, with dual $\mathbb{X}_j = (X_j,R_j)$. Then the dual of $\prod \mathbb{A}_j$ is 
\[ \mathbb{X} = ( \beta(\Sigma X_j), \overline{\Sigma R_j}) \] where $\Sigma$ is cardinal sum and $\beta(\sum X_j)$ the Stone-\v{C}ech compactification of $\sum X_j$ endowed with the discrete topology.
\end{lemma}
\begin{proof}
The theorem is well known for the underlying Stone spaces. It remains to prove that the accessibility relation - let us denote it by $R^\beta$ - in the dual of $\prod \mathbb{A}_j$ is the closure of $\sum R_j$. The procedure is similar to the one used in the proof of Proposition \ref{prop_dual_modalisation}.
\end{proof}

\begin{lemma}\label{lem2.26} For any family $(\mathbb{A}_j \mid j \in J)$ of subordination algebras, there is a unique onto morphism in \textbf{sCM} $f : (\prod \mathbb{A}_j)^\delta \longrightarrow \prod \mathbb{A}_j^\delta$ such that $\id \circ f = p$, where $\id$ is the identity map $\prod \mathbb{A}_j \longrightarrow (\prod \mathbb{A_j})^\delta$ and $p$ is the canonical weak embedding $\prod \mathbb{A}_j \longrightarrow \prod \mathbb{A}_j^\delta$.
\end{lemma}
\begin{proof}
Each projection $p_j : \prod \mathbb{A}_j \longrightarrow \mathbb{A}_j$ lifts to a morphism $p_j^\delta$ in \textbf{sCM} by \ref{product}, $(\prod \mathbb{A}_j)^\delta \longrightarrow \mathbb{A}_j^\delta$, and the product of the $p_j^\delta$'s is the required morphism $f$. If the dual of $\mathbb{A}_j$ is $\mathbb{X}_j$, then by \ref{prop2.25}, we may write $(\prod \mathbb{A}_j)^\delta = \mathcal{P}(\beta(\sum X_j))$, $\prod \mathbb{A}_j^\delta = \mathcal{P}(\sum X_j)$ and $f$ is the map $E \longmapsto E \cap \sum X_j$ which is clearly onto.
\end{proof}

\begin{definition} The map defined in Lemma \ref{lem2.26} is the canonical epimorphism $(\prod \mathbb{A}_j)^\delta \longrightarrow \prod \mathbb{A}_j^\delta$. Its restriction $f^m$ to $(\prod \mathbb{A}_j)^m$ clearly takes its values in $\prod \mathbb{A}_j^m$. We call it the canonical morphism (in \textbf{MA}) $(\prod \mathbb{A}_j)^m \longrightarrow \prod \mathbb{A}_j^m$. We say that $(\mathbb{A}_j \mid j \in J)$ is a \textbf{good} family if the canonical map $f^m$ is an embedding.
\end{definition}

\begin{proposition}\label{goodfamproduct} Let $(\mathbb{A}_j \mid j \in J)$ be a family of subordination algebras. Then $\mathbb{A}_j \models \varphi$ for all $j$ implies $\prod \mathbb{A}_j \models \varphi$ for all modal formulas $\varphi$ if and only if $(\mathbb{A}_j \mid j \in J)$ is a good family.
\end{proposition}
\begin{proof}
The elements of $(\prod A_j)^m$ are of the form $\varphi(a_1,...,a_n)$ for some modal formula $\varphi$ and elements $a_1,...,a_n$ of $\prod A_j$. Since $f^m$ is a modal algebra homomorphism, then 
\[ f^m(\varphi(a_1,...,a_n)) = (\varphi(a_{1j},...,a_{nj})\mid j \in J). \] And $f^m$ is an embedding if and only if its $1$-kernel is reduced to $\lbrace 1 \rbrace$, that is, $\varphi(a_1,...,a_n) = 1$ if and only if $\varphi(a_{1j},...,a_{nj}) = 1$ for all $j$, for all $\varphi$ and all $a_1,...,a_n$. Now $\varphi(a_1,...,a_n) = 1$ means $\prod \mathbb{A}_j \models_v \varphi$ for the valuation $v$ sending the variable $p_i$ to $a_i$, $i = 1,...,n$, and the proposition is proved.
\end{proof}

\section{Correspondence theory}\label{Section3}

We now turn to problems of correspondence theory. The context of subordination algebras is suitable to study the interconnection between three languages: the (bi)modal one, the subordination one and the accessibility one. Hence, we have several kinds of correspondence that may be studied.

First of all, now that subordination algebras have been established as models for modal logics, we have the classical correspondence aspect which is concerned with the translation of modal (or bimodal) equations on a subordination algebra into first-order properties of the accessibility relation on its dual. For instance, one can show that, for a subordination algebra $\mathbb{B}$, we
 have that $ \mathbb{B} \models \square p \rightarrow \square \square p$ if and only if its associated accessibility relation is transitive. 
 
 Then, a second kind of correspondence may be studied: the one concerned with the translation of first-order properties of the subordination language into first order properties of the accessibility relation. For instance, a subordination algebra satisfies axiom S8, that is $a \prec b$ implies $a \prec c \prec b$ for some $c$, if and only if its associated accessibility relation is transitive.  We will not consider this theory here and redirect to \cite{Balbiani} and \cite{SantoliThese} for results in this direction.
 
Finally, a third kind of correspondence (which could be glimpsed in \cite{Vakar}) naturally arises from the previous ones. Indeed, one can look for translate modal or bimodal equations into first-order properties in the language of the subordination algebras (i.e. using the Boolean connectives and the subordination $\prec$). Directly following from our previous examples, we have for instance that $\mathbb{B} \models \square p \rightarrow \square \square p$ if and only if $\mathbb{B}$ satisfies axiom S8.

We first begin with some specific examples of translation, to get the flavour of more general results. Our examples will be given in the realm of bimodal formulas. Recall from \ref{notationscheme} that $\varphi(\overline{\psi})$ is the scheme associated to the formula $\varphi(\overline{p})$ and that $\mathbb{B} \models \varphi(\overline{\psi})$ means $\mathbb{B} \models \varphi(\overline{\psi})$ for every tuple $\overline{\psi}$. We begin by an idyllic example.

\begin{lemma}\label{Impexem}For a subordination algebra $\mathbb{B}$ with dual $\mathbb{X}$ and for any $k,l,m,n\in \N$, the following are equivalent: 
\begin{enumerate}
\item $\mathbb{B} \models \lozenge^k \square^l p \rightarrow \square^m \lozenge^n p$,
\item $\mathbb{B} \models \blacklozenge^k \lozenge^m p \rightarrow \lozenge^l \blacklozenge^n p$,
\item $\mathbb{B} \models \lozenge^k \square^l \psi \rightarrow \square^m \lozenge^n \psi$,
\item $\mathbb{B} \models \blacklozenge^k \lozenge^m \psi \rightarrow \lozenge^l \blacklozenge^n \psi$,
\item $\mathbb{X} \models (x \mathrel{R^k}y \text{ and } x \mathrel{R}^mz) \rightarrow (\exists u)(y \mathrel{R}^lu \text{ and } z \mathrel{R}^n u)$
\item $\mathbb{B} \models (\neg a \prec^l \neg b \text{ and } a \perp^n c) \rightarrow (\exists d)( b\prec^k d \text{ and } c \perp^m  d)$, 
\end{enumerate}
where $a \perp b$ is a shortcut for $ a \prec \neg b$.
\end{lemma}
\begin{proof}
Of course, $3. \Rightarrow 1.$ and $4. \Rightarrow 2.$.

We know consider the following sequence of equivalences, in which $A,B,C$ and $D$ are clopen:  
\begin{enumerate}
\item $R^{-k}[ R^{-l}[ A^c]^c] \subseteq R^{-m}[ R^{-n}[A]^c]^c$,
\item $R^{-k}[ R^{-l}[ A^c]^c] \cap R^{-m}[ R^{-n}[A]^c] = \emptyset$,
\item For all $B$ and $C$ such that $B  \subseteq  R^{-l}[ A^c]^c$ and  $ C \subseteq  R^{-n}[A]^c$, we have that $R^{-k}[B] \cap R^{-m}[C] = \emptyset$,
\item For all $B$ and $C$ such that $B  \subseteq R^{-l}[ A^c]^c$ and  $ C \subseteq  R^{-n}[A]^c$, there exists $ D$ such that $R^{-k}[B] \subseteq D$ and $R^{-m}[C] \subseteq  D^c$,
\item For all $B$ and $C$, $ B^c  \supseteq R^{-l}[ A^c]$ and $(C \times A) \cap R^n = \emptyset$ imply there exists $D$ such that $R^{-k}[B] \subseteq D$ and $(D \times C) \cap R^m = \emptyset$.
\end{enumerate}
This sequence shows that $1.$ (i.e. a)) and $6.$ (i.e. e)) are equivalent.

We now prove $6. \Rightarrow 5.$ Suppose $x\mathrel{R}^k y$ and $x \mathrel{R}^mz$ while for no $u$, $y\mathrel{R}^lu$ and $z \mathrel{R}^n u$. Then $R^l(y) \cap R^n(z) = \emptyset$. So there is $A$ with $R^l(y) \subseteq  A$ and $R^n(z) \subseteq  \neg A$. And there is $B \ni y$ and $C \ni z$ such that $R^l[B] \subseteq  A$ and $R^n[C] \subseteq A^c$, in other words with $ R^{-l}[ A^c] \subseteq  B^c$ and $C \subseteq  R^{-n}[ A]^c$. By d) (which is equivalent to $6.$)), there is $D$ with $R^{-k}[B] \subseteq D$ and $R^{-m}[C] \subseteq \neg D$. Since $x \mathrel{R}^k y$, we have $x \in D$ and since $y \mathrel{R}^mz$, we have $x \in  D^c$, which is impossible.

Finally, by modal logic, for the accessibility condition $G_{klmn}$ (as it is denoted in \cite{Chellas}) $5.$, we have the following equivalences
\begin{align*}
&\mathbb{X} \models G_{klmn} \\
\iff &\mathbb{X} \models \lozenge^k \square ^l p \rightarrow \square ^m \lozenge^n p \\
\iff & \mathbb{B}^\delta \models \lozenge^k \square ^l p \rightarrow \square ^m \lozenge^n p  \\
\iff & \mathbb{B}^\delta \overset{\text{\small{schm}}}{\models} \lozenge^k \square^l \psi \rightarrow \square^m \lozenge^n \psi.
\end{align*}
 In other words, we have that $5.$ implies $3.$.

We have proved $3. \Rightarrow 1. \Rightarrow 6. \Rightarrow 5. \Rightarrow 3.$. One proves $4. \Rightarrow 2.\Rightarrow 6. \Rightarrow 5. \Rightarrow 4.$ in a similar way.
\end{proof}

In the example, the characterisation of modal formulas in term of the accessibility relation is exactly the same as in the purely modal case. We give a two variables example of this phenomenon (without proof since this example, as well as in \ref{Impexem}, is taken into account in Theorems \ref{Sahlq} and \ref{Sahlq2}).

\begin{example}\label{Impexem2}For a subordination algebra $\mathbb{B}$ with dual $\mathbb{X}=(X,R)$ the following are equivalent:
\begin{enumerate}
\item $\mathbb{B} \models \square(\square p \rightarrow q) \vee \square (\square q \rightarrow p)$,
\item $(X,R) \models (x \mathrel{R} y \text{ and } x \mathrel{R} z) \rightarrow (y \mathrel{R} z \text{ or } z \mathrel{R} y)$,
\item $\mathbb{B} \models( a \perp b \text{ and } b \perp a) \rightarrow ((\exists c)(a \prec c \text{ and } b \perp c))$.
\end{enumerate}
\end{example}

We now give an analogue of Sahlqvist theorem, that is, give a set of modal formulas that are first-order expressible in a uniform way. In particular, those formulas are subordination canonical, in a natural sense given in \ref{defcanonic}. The obtained set of Sahlqvist formulas for subordination algebras is definitely smaller than the set of Sahlqvist formulas for modal algebras. This is justified by the fact that there exists (see \ref{counterexemp}) Sahlqvist formulas which are not subordination canonical.

\begin{definition} A bimodal formula $\varphi$ is \textbf{closed} (resp. \textbf{open}) if it is obtained from constants $\top$, $\bot$, propositional variables and their negations, by applying $\vee$, $\wedge$, $\lozenge$ and $\blacklozenge$ (resp. $\vee$, $\wedge$, $\square$ and $\blacksquare$).

A bimodal formula $\varphi$ is \textbf{positive} (resp. \textbf{negative}) if it is obtained from constants $\top$, $\bot$ and propositional variables (resp. and negations of propositional variables) by applying $\wedge$, $\vee$, $\lozenge$, $\square$, $ \blacklozenge$ and $\blacksquare$.

A bimodal formula $\varphi$ is \textbf{s-positive} (resp. \textbf{s-negative}) if it is obtained from closed positive formulas (resp. open negative formulas) by applying $\vee$, $\wedge$, $\square$ and $\blacksquare$ (resp. $\vee$, $\wedge$, $\lozenge$ and $\blacklozenge$).

A bimodal formula $\varphi$ is \textbf{g-closed} (resp. \textbf{g-open}) (g for generalised) if it is obtained from closed (resp. open) formulas by applying  $\vee$, $\wedge$, $\square$ and $\blacksquare$ (resp. $\vee$, $\wedge$, $\lozenge$ and $\blacklozenge$).
\end{definition}

To obtain the analogue of Sahlqvist result, we need two more ingredients.

\begin{definition}\label{DefSahl} A \textbf{strongly positive} bimodal formula is conjunction of formulas of the form 
\[ \square^{\langle \mu \rangle} p := \square^{\mu_1}\blacksquare^{\mu_2}\cdots \square^{\mu_k} p , \] where $p \in \Var$ and $\mu \in \N^k$ for some $k \in \N$.

A \textbf{s-untied} bimodal formula is a formula obtained from strongly positive and s-negative formulas by applying only $\wedge$, $\lozenge$ and $\blacklozenge$.

Finally, a formula $\varphi$ is said to be \textbf{s-Sahlqvist} if of the form $\varphi = \square^{\langle \mu \rangle}(\varphi_1 \rightarrow \varphi_2)$ where $\varphi_1$ is s-untied and $\varphi_2$ s-positive. By definition any s-Sahlqvist formula $\varphi$ is a Sahlqvist formula and by Sahlqvist's theorem (\cite{SahlqvistOrigin}, \cite{Sahlqvist} and adapted for bimodal formulas (among others) in \cite{ExtSahlq}), there is a first order formula $f(\varphi)$ in the language of a binary relation such that for any bimodal algebra $\mathbb{B}$ with dual $\mathbb{X}$, $\mathbb{B} \models \varphi$ if and only if $\mathbb{X} \models f(\varphi)$.
\end{definition}

\begin{theorem}\label{Sahlq} Let $\varphi$ be an s-Sahlqvist bimodal formula and let $f(\varphi)$ be its associated first-order formula as defined in Definition \ref{DefSahl}. Then for any subordination algebra $\mathbb{B}$ with dual $\mathbb{X}$, we have 
\[ \mathbb{B} \models \varphi \text{ if and only if } \mathbb{X} \models f(\varphi). \]
\end{theorem}
\begin{proof}
We prove Sahlqvist theorem in the generalised context of subordination algebras simply by following the topological proof of Sambin and Vaccaro in \cite{Sahlqvist}. In almost all places only the closedness of the accessibility relation is needed. The only place where the extra assumption that $R^{-1}[O]$ is open when $O$ is open is necessary is in the intersection lemma (see below). This explains our definition of s-Sahlqvist formulas. The intersection lemma we then use is the following.

\paragraph{Intersection lemma}Let $\varphi(p_1,...,p_k)$ be an s-positive bimodal formula and $\mathbb{X}=(X,R)$ a subordination space. For every $A \subseteq X$ and for every $C_1,...,C_{k-1}$ closed sets of $X$
\[ \varphi(C_1,...,\cl(A),...,C_{k-1}) = \bigcap \lbrace \varphi(C_1,...,O,...,C_{k-1}) \mid  A \subseteq O \in \Of(X) \rbrace,\] where $\cl(A)$ denotes the topological closure of $A$.

\paragraph{Proof of the lemma}The proof is done by induction on the complexity of $\varphi$. We note that, $\varphi$ being s-positive,  $\varphi \equiv \lozenge \psi$ or $\blacklozenge \psi$ implies that $\psi$ does not contain any white or black boxes. Hence, $\psi$ is a closed formula, and it follows that \[( \psi(C_1,...,O,...,C_{k-1}) \mid A \subseteq O \in \Of(X) ) \] is a filtered family of closed sets. This allows us to use Esakia's Lemma as in \cite{Sahlqvist} and conclude the proof.
\end{proof}

\begin{definition}\label{defcanonic} A bimodal formula $\varphi$ is \textbf{s-canonical} if $\mathbb{B} \models \varphi$ implies $\mathbb{B}^\delta \models \varphi$ for any subordination algebra $\mathbb{B}$. it is said to be \textbf{scheme-extensible} if $\mathbb{B} \models \varphi(\overline{p})$ (we write $\varphi(\overline{p})$ to indicate that the variables of $\varphi$ are among the tuple $\overline{p}$) implies $\mathbb{B} \models \varphi(\overline{\psi})$ for all $\overline{\psi}$. Clearly, the latter is equivalent to $\mathbb{B} \models \varphi$ implies $\mathbb{B}^m \models \varphi$. Hence, being s-canonical implies being scheme-extensible (since $\mathbb{B}^m$ is a subalgebra of $\mathbb{B}^\delta$).
\end{definition}

\begin{corollary} Any s-Sahlqvist bimodal formula is s-canonical and therefore is scheme-extensible.
\end{corollary}
\begin{proof}
Let $\varphi$ be an s-Sahlqvist formula. Then $\varphi$ is a Sahlqvist formula and $(X,R) \models \varphi$ if and only if $(X,R) \models f(\varphi)$ for any Kripke frame $(X,R)$. Therefore, $\mathbb{B} \models \varphi$ if and only if $\mathbb{X} \models f(\varphi)$ (if $\mathbb{X} = (X,R,\tau)$ is the dual of $\mathbb{B})$ by Theorem \ref{Sahlq}. The latter being equivalent to $(X,R) \models \varphi$ if and only if $\mathbb{B}^\delta \models \varphi$.
\end{proof}

\begin{example}\label{counterexemp} Let us have a look at the formula 
\[ \varphi \equiv p \rightarrow \lozenge \square p, \]
already examined in Example \ref{exampl}. It is a Sahlqvist formula, but not an s-Sahlqvist formula. On modal algebras, it is equivalent to the formula 
\[ f(\varphi) \equiv (\forall x)(\exists y)(x \mathrel{R} y \text{ and } R(y) \subseteq \lbrace x \rbrace). \] This fact is no longer true on subordination algebras: the subordination space of Example \ref{exampl} satisfies $\varphi$ but not $f(\varphi)$. 

Finally, the formula $\varphi$ is not scheme-extensible, as shown in \ref{exampl}, hence $\varphi$ is an example of canonical formula which is not s-canonical.
\end{example}

We now study our second kind of correspondence theory, namely the translation of a bimodal formula into the subordination algebra language. Examples of such translations have already been given in Examples \ref{Impexem} and \ref{Impexem2}. As promised, we generalise these examples in the next results.

\begin{definition}\label{defdef} A bimodal formula $\varphi = \varphi(\overline{p})$ is said to be \textbf{s-definable} (resp. \textbf{$\leq$-definable\text{;} $\geq$-definable}) if there is an effectively produced first order formula $\xi = \xi(\varphi) = \xi(\overline{p})$ (resp. $\xi_\leq=\xi_\leq(\varphi) = \xi_\leq(\overline{p},q)$ and $\xi_\geq=\xi_\geq(\varphi) = \xi_\geq(\overline{p},q)$) such that for any subordination algebra $\mathbb{B}$ and any valuation $v : \Var \longrightarrow \mathbb{B}$, one has: 
\begin{enumerate}
\item $\mathbb{B} \models_v \varphi(\overline{p})$ if and only if $\mathbb{B} \models_v \xi(\overline{p})$,
\item $\mathbb{B} \models_v \varphi(\overline{p})\rightarrow q^\pm$ if and only if $\mathbb{B} \models_v \xi_\leq(\overline{p};q^\pm)$ (where $q^\pm$ is a shorthand for $q$ or $\neg q$),
\item $\mathbb{B} \models_v q^\pm \rightarrow \varphi(\overline{p})$ if and only if $\mathbb{B} \models_v \xi_\geq(\overline{p},q^\pm)$.
\end{enumerate}
Clearly, if $\varphi$ is $\geq$-definable (resp. $\leq$-definable), then $\phi$ (resp. $\neg \varphi$) is s-definable. Also $\varphi$ is $\geq$-definable if and only if $\neg \varphi$ is $\leq$-definable.
\end{definition}

\begin{theorem} \label{prelasttheorem}If $\varphi$ is an open or a closed formula, then both $\varphi$ and $\neg \varphi$ are both $\leq$ and $\geq$-definable.
\end{theorem}
\begin{proof}
We begin by the following general remark, that will help to facilitate computation. We may assume that our working subordination algebra is $\mathbb{B} = \Of(\mathbb{X})$ where $\mathbb{X}=(X,R)$ is the dual of $\mathbb{B}=(B,\prec)$. Under a valuation $v$, variables $p$ and their negations $\neg p$ are therefore clopen subsets of $X$ and more generally, formulas are subsets of $X$. Also, $\square \varphi =  R^{-1}[ \neg \varphi]^c$, $\blacksquare \varphi =  R[\neg \varphi]^c$, $ \lozenge \varphi = R^{-1}[ \varphi]$ and $\blacklozenge \varphi = R[\varphi]$. On the subordination side, remember that $\varphi \prec \psi$ is equivalent to each of the following conditions: $R^{-1}[\varphi] \subseteq \psi$ ; $\neg \psi \times \phi \cap R = \emptyset$ and $\varphi \subseteq  R[\neg \psi]^c$. Hence, each of these expressions, when restricted to clopen subsets $p^\pm$ (see the second point of Definition \ref{defdef}) of $X$, corresponds to an atomic formula in the first order language of subordination algebra. Finally, we make use of the following topological remarks, in which $A \subseteq X$, $O$ is an open subset and $F$ closed subset of $X$: 
\begin{enumerate}
\item $O \subseteq A$ if and only if for all variables $p$, $p \subseteq O$ implies $p \subseteq A$,
\item $A \subseteq F$ if and only if for all variables $p$, $F \subseteq p$ implies $A \subseteq p$,
\item $F \subseteq O$ if and only if for some variable $p$, $F \subseteq p \subseteq O$,
\item $R^{-1}[F] \subseteq O$ if and only if for some variables $p,q$, one has $F \subseteq p$, $q \subseteq O$ and $R^{-1}[p] \subseteq q$.
\end{enumerate}

We are ready for the proof, that is done by induction on the complexity of $\varphi$. We only consider the case where $\varphi$ is open since $\varphi$ is closed if and only if $\neg \varphi$ is open.

If $\varphi$ is a constant, a variable or the negation of a variable, the result is clear by our beginning remark.

Suppose now $\varphi \equiv \theta \vee \psi$. Then $\varphi \rightarrow q^\pm$ is equivalent to $(\theta \rightarrow q^\pm) \wedge (\psi \rightarrow q^\pm)$ and the result follows by the induction hypothesis. Also, $q^\pm \rightarrow \theta \vee \psi$ is equivalent to 
\[(\exists r,s)((r\rightarrow \theta) \wedge (s \rightarrow \psi) \wedge (q^\pm \rightarrow r \vee s)). \] Here again, we use induction to conclude.

If $\varphi \equiv \theta \wedge \psi$, then $q^\pm \rightarrow \varphi$ is equivalent to $(q^\pm \rightarrow \theta) \wedge (q^\pm \rightarrow \psi)$ while $\varphi\rightarrow q^\pm$ is equivalent to 
\[ (\forall r)(((r\rightarrow \theta) \wedge (r \rightarrow \psi))\rightarrow (r \rightarrow q^\pm)). \]

Finally, we consider the case $\varphi \equiv \square \psi$. Then $\varphi \rightarrow q^\pm$ is equivalent to 
\[ (\forall r)(r \subseteq R^{-1}[\neg \psi]^c \rightarrow r \subseteq q^\pm), \] which is in turn equivalent to 
\[ (\forall r)((\exists s)((\neg \psi \subseteq s \text{ and } R^{-1}[s] \subseteq \neg r) \rightarrow (r \rightarrow q^\pm)). \] And $q^\pm \rightarrow \square \psi$ is equivalent to $q^\pm \subseteq R^{-1}[\neg \psi]^c$, that is $R^{-1}[\neg \psi] \rightarrow q^\mp$, which is equivalent to 
\[ (\exists r)(\neg \psi \subseteq r \text{ and } R^{-1}[r] \subseteq q^\mp) . \]
\end{proof}

We now need an analogue of the intersection lemma.

\begin{lemma}\label{lastlem} If $\varphi$ is positive open and $\overline{\psi}$ is closed, then 
\[ \varphi(\overline{\psi}) = \bigcap \lbrace\varphi(\overline{p}) \mid \overline{p}\text{ clopen  and } \overline{p} \geq \overline{\psi} \rbrace. \]
\end{lemma}
\begin{proof}
Since $\varphi$ is positive, we have $\varphi(\overline{\psi}) \subseteq \bigcap \lbrace \varphi(\overline{p}) \mid \overline{p} \geq \overline{\psi}\rbrace$ and we prove the opposite inclusion $\supseteq$ by induction on the complexity of $\varphi$.

This is clear when $\varphi$ is a variable, because $\overline{\psi}$ is closed. 

Consider the case $\varphi \equiv \xi \vee \theta$. If $x \in \varphi(\overline{p})$ for all $\overline{p} \geq \overline{\psi}$ but $x \nin \varphi(\overline{\psi})$, then $x \nin \xi(\overline{\psi})$ and $x \nin \theta(\overline{\psi})$. By induction, there $\overline{p} \geq \overline{\psi}$ with $x\nin \xi(\overline{p})$ and $\overline{q} \geq \overline{\psi}$ with $x \nin \theta(\overline{q})$. Then $\overline{p } \cap \overline{q}$ is clopen and $\overline{\psi} \leq \overline{p } \cap \overline{q}$ so that 
\[ x \in \varphi(\overline{p } \cap \overline{q}) = \xi(\overline{p } \cap \overline{q}) \cap \theta(\overline{p } \cap \overline{q}) \subseteq \xi(\overline{p}) \cap \theta(\overline{q}), \] a contradiction.

If $\varphi \equiv \xi \wedge \theta$, then \[ \varphi(\overline{\psi}) = \xi(\overline{\psi}) \cap \theta(\overline{\psi}) = \bigcap \lbrace \xi(\overline{p}) \mid \overline{p} \geq \overline{\psi}\rbrace \cap \bigcap \lbrace \theta(\overline{p}) \mid \overline{p} \geq \overline{\psi}\rbrace = \bigcap \lbrace \varphi(\overline{r}) \mid \overline{r} \geq \overline{\psi}\rbrace .\]

Finally, suppose $\varphi \equiv \square \theta$. Then, 
\begin{align*}
\varphi(\overline{\psi}) &=  R^{-1}[\neg \theta(\overline{\psi})]^c\\
&=  R^{-1}[\bigcup\lbrace \neg \theta(\overline{p}) \mid \overline{p} \geq \overline{\psi} \rbrace]^{c}\\
&=\bigcap \lbrace  R^{-1}[\neg \theta(\overline{p})]^c \mid \overline{p} \geq  \overline{\psi} \rbrace\\
&= \bigcap \lbrace \varphi(\overline{p}) \mid \overline{p} \geq \overline{\psi} \rbrace
\end{align*} as required.
\end{proof}

\begin{theorem}\label{Sahlq2} If $\xi$ is a g-closed formula, then $\xi$ is $\geq$-definable, hence s-definable.
\end{theorem}
\begin{proof}
If $\xi$ is g-closed, there is a positive open formula $\varphi$ and a  tuple of closed formulas $\overline{\psi}$ such that $\xi = \varphi(\overline{\psi})$. Then, the formula $q^\pm \rightarrow \varphi(\overline{\psi})$ is equivalent, by Lemma \ref{lastlem}, to 
\[ \forall \overline{p} \geq \overline{\psi}, q^\pm \rightarrow \varphi(\overline{p}). \] And both formulas $\overline{p} \geq \overline{\psi}$ and $q^\pm \rightarrow \varphi(\overline{p})$ are s-definable by Theorem \ref{prelasttheorem}.
\end{proof}

As announced, we now compare the three modal languages (white, black and bicolour) one with another. The comparison is first done semantically by establishing analogues of Birkhoff's characterisation of varieties for each modal language. Specific examples are then derived.

In universal algebra, Birkhoff theorem is twofold. First, a characterisation of those sets of identities which are true in a class of algebras in term of a provability system. And then, a characterisation of those classes of algebras that satisfy some set of identities in terms of semantic constructs.

In modal algebra, where identities may be assimilated to formulas, our provability system always gives a normal modal logic as set of theorems and as discussed in Definition \ref{defval}, this is not always the case for the logic of a class $\mathcal{K}$ of subordination algebras. So our first result will be a criterion to ensure that $\mbox{Log}(\mathcal{K})$ is a normal modal logic and, then, give a characterisation of those classes of subordination algebras that satisfy some normal modal logic.

\begin{proposition} If $\mathcal{K}$ is a class of subordination algebras, then $L = \mbox{Log}(\mathcal{K})$ is a normal modal logic if and only if $\mathbb{B} \in \mathcal{K}$ implies $\mathbb{B}^m \in \mathcal{K}$.
\end{proposition}
\begin{proof}
Suppose $L = \mbox{Log}(\mathcal{K})$ is a normal modal logic and $\mathbb{B} \in \mathcal{K}$. Since $L$ may be axiomatized by schemes, this follows directly from Proposition \ref{modal1}. 

Suppose now $\mathbb{B} \in \K$ implies $\mathbb{B}^m \in \K$. We have to prove that $L$ is closed under substitution, that is $\varphi \in L$ implies $\varphi(\overline{\psi}) \in L$. Let $\mathbb{B} \in \K$. Then, $\mathbb{B}^m \in \K$ and so $\mathbb{B}^m \models \varphi$, whence $\mathbb{B}^m \models \varphi(\overline{\psi})$ as $\mathbb{B}^m$ is a modal algebra, and it follows that $\mathbb{B} \models \varphi(\overline{\psi})$ as proved in \ref{modal1}.
\end{proof}

\begin{theorem}\label{whiteBirk} Let $\K$ be a class of subordination algebras. Then the following are equivalent: 
\begin{enumerate}
\item $\K = \mo(L)$ for some modal normal logic $L$,
\item $\K$ is definable by schemes of modal formulas,
\item $\K$ is closed under subalgebras and morphic images (in {\bf{Sub}}), products of good families and modalisations, and reflects modalisation, that is $\mathbb{B}^m \in \K$ implies $\mathbb{B} \in \K$, 
\item $\K$ is closed under subalgebras and morphic images (in {\bf{Sub}}), and for any family $(\mathbb{B}_i \mid i \in I)$, one has 
\[ \prod \mathbb{B}_i^m \in \K \text{ if and only if } \forall i \in I, \ \mathbb{B}_i \in K. \]
\end{enumerate}
\end{theorem}
\begin{proof}
The equivalence $1. \Leftrightarrow 2.$ is clear. Both implications $2. \Rightarrow 3.$ and $2. \Rightarrow 4.$ follows for \ref{modal1}, \ref{quotient}, \ref{subalgebra} and \ref{goodfamproduct}.

Let us prove $3. \Rightarrow 2.$ (one proves $4. \Rightarrow 2.$ in a similar way). Let $\mathcal{M} = \lbrace \mathbb{B} \in \textbf{MA} \mid \mathbb{B} \in \K \rbrace$. Then $\mathcal{M}$ is a class of modal algebras closed under H,S and P and is therefore an equational class by Birkhoff classical theorem. Let $L$ be an axiomatisation of $\mathcal{M}$ by schemes. All we have to prove is $\K = \mo(L)$.

If $\mathbb{B} \in \mo(L)$, then $\mathbb{B}^m \in \mo(L)$ by Proposition \ref{modal1}, so that $\mathbb{B}^m \in \mathcal{M} \subseteq \K$. Since $\K$ reflects modalisation, it follows that $\mathbb{B} \in \K$. Conversely, if $\mathbb{B} \in \K$, then $\mathbb{B}^m \in \K$ and, being a modal algebra, $\mathbb{B}^m \in \mathcal{M}$. Hence, $\mathbb{B}^m \in \mo(L)$. It follows from \ref{modal1} that $\mathbb{B} \in \mo(L)$.
\end{proof}

Of course, there is a black and a bimodal version of this theorem. We only present the bimodal version. A family $(B_i \mid i \in I)$ of subordination algebras is said to be \textbf{s-good} if the canonical morphism $f^{bim}$ (in \textbf{sSub}) (the restriction $(\prod \mathbb{B}_i)^{bim} \longrightarrow \prod \mathbb{B}_i^{bim}$ of the canonical epimorphism $(\prod \mathbb{B}_i)^\delta \longrightarrow \prod \mathbb{B}_i^\delta$) is an embedding. Of course, this is stronger than being good. 

\begin{theorem}\label{biBirk} Let $\K$ be a class of subordination algebras. Then, the following are equivalent: 
\begin{enumerate}
\item $\K = \mo(L)$ for some bimodal tense logic $L$,
\item $\K$ is closed under strong subalgebras, strong morphic images, product of s-good families, bimodalisations and reflects bimodalisation.
\end{enumerate}
\end{theorem}

This leaves open the non-scheme versions of the two theorems.

\begin{problem} \begin{enumerate}
\item Characterise the sets of formulas of the form $\Log(\K)$ where $\K$ is a class of subordination algebras in term of provability - and give the associated completeness theorem.
\item Characterise semantically the equational classes of subordinations algebras, that is, the classes $\mo(L)$ where $L$ is an arbitrary set of modal formulas (not necessarily closed under substitution).
\end{enumerate}
\end{problem}

A fourth kind of correspondence can be realised within the realm of modal formulas, if we remember that, for unimodal formulas, three different languages may be adopted: the white language, the black one and the bicolour (bimodal) one. An example of this phenomenon is given in example \ref{Impexem}: the bicolour formula $\lozenge p \rightarrow \blacklozenge p$ is equivalent to the white formula $\lozenge \square p \rightarrow p$, and to the black formula $\blacklozenge \blacksquare p \rightarrow p$ (all are equivalent to the symmetry of $R$). At the theoretical level, everything is settled by the (white, black and bicolour) Birkhoff theorems (\ref{whiteBirk} and \ref{biBirk}), and we just give here some examples and counterexamples of correspondences  between these three languages.

\begin{example} We know by Example \ref{Impexem} that the accessibility condition $G_{klmn}$
\[ (x\R^k y \text{ and } x \R^m z) \rightarrow (\exists u)(y \R^l u \text{ and } z \R^n u) \] is equivalent to the bicolor formula $\blacklozenge^k \lozenge^m p \rightarrow \lozenge^l \blacklozenge ^n p$ and even to a white formula ($\lozenge^k \square^l p \rightarrow \square^m \lozenge ^np$). It is equivalent to a black modal formula only in the cases where $k$ or $m=0$, and $l$ or $n=0$.

Since the condition is symmetrical under the change $k$ with $m$ and $l$ with $n$, it suffices to examine the cases $m=l=0$ and $m=n=0$.

If $m=l=0$, then the condition is trivially equivalent to $\blacklozenge^k p \rightarrow \blacklozenge^n p$.

If $m=n=0$, the condition is equivalent to $\blacklozenge^k \blacksquare^l p \rightarrow p$.

We  now prove that in the other cases, there is no black equivalent formula.

Suppose first $k,m > 0$. We give an example of subordination space $\mathbb{X}=(X,R)$ such that $\mathbb{X} \models G_{k0l0}$ but which admits a \textbf{$\blacklozenge$SubS} morphic image $\mathbb{Y}$ such that $\mathbb{Y} \not\models G_{k0m0}$. It suffices to take $\mathbb{X} = \lbrace a,b,c,d \rbrace$, $R = \lbrace (a,b), (b,b), (c,d), (d,d) \rbrace$, $\Theta$ the congruence which relates $a$ with $c$ and $\mathbb{Y} = \mathbb{X}/\Theta$.

Suppose now $l,n > 0$. We give an example of a subordination space $\mathbb{X}=(X,R)$ such that $\mathbb{X} \models G_{0l0n}$ but which admits a \textbf{$\blacklozenge$SubS} subobject $\mathbb{Y} \not\models G_{0l0n}$. We take $X= \lbrace x_0,x_1,\cdots,x_{l},y_1,\cdots y_{n-1} \rbrace$, $R = \lbrace (x_i,x_{i+1}) \mid i = 0,...,l-1 \rbrace \cup \lbrace (x_i,x_i) \mid i = 1,...,l \rbrace \cup \lbrace (x_0,y_1),(y_{n-1},x_{l} \rbrace \cup \lbrace (y_i,y_{i+1} \mid i=1,...,n-1\rbrace$ and $\mathbb{Y}$ is the subobject given by $\lbrace x_0 \rbrace$.
\end{example}

\begin{example} The axiom $\lozenge \blacklozenge \lozenge p \rightarrow \lozenge p$ is not expressible by a unicolour axiom.
\end{example}
\begin{proof}
It is not difficult to see that the mentioned axiom correspond to the first order property 
\begin{equation}\label{eqexem1} x \mathrel{R} y, z \mathrel{R} y, z \mathrel{R} u \rightarrow x \mathrel{R} u. \end{equation}
Let $\mathbb{X} = (X,R)$ be the subordination space with $X = \lbrace a,b,c,d,e \rbrace$ and $R = \lbrace (a,b),(c,d),(c,e) \rbrace$. Then, $\mathbb{X} \models \eqref{eqexem1}$ (vacuously). Now, the equivalence $\Theta$ generated by $\lbrace (b,d) \rbrace$ is a congruence in $\textbf{SubS}$. The quotient $\mathbb{X}/\Theta$ is $(X/\Theta,R/\Theta)$ where $X/\Theta = \lbrace a^\Theta, b^\Theta = d^\Theta, c^\Theta, e^\Theta \rbrace$ and $R/\Theta = \lbrace (a^\Theta,b^\Theta), (c^\Theta,b^\Theta), (c^\Theta,e^\Theta)\rbrace$ and clearly $\mathbb{X}/\Theta \not\models \eqref{eqexem1}$. This shows that axiom \eqref{eqexem1} is not expressible in the white language. One prove in a similar way (consider $(X,R^\partial)$) that axiom \eqref{eqexem1} is not expressible in the black language.
\end{proof}

We end by two examples of formulas in the subordination language which are not modally definable.

\begin{example} The connectedness axiom 
\begin{description}
\item[{\normalfont (Con)}] $a \prec a$ implies $a = 0$ or $a =1$ 
\end{description} is not modally definable, even in the bicolour language.
\end{example}
\begin{proof}
Let $\mathbb{B}$ be a non-trivial subordination algebra satisfying (Con) (see \cite[Theorem 2.6.1]{Vakar} for examples of such algebras). Then $\mathbb{B} \times \mathbb{B}$ does not satisfies (Con) since $(1,0) \neq (0,0),(1,1)$ but $(1,0) \prec (1,0)$. It follows that (Con) is not modally definable. 
\end{proof}

\begin{example} The extensionality axiom 

\begin{description}
\item[{\normalfont (S5)}] $ a \neq 0$ implies $b \prec a$ for some $b \neq 0$,
\end{description}
is not modally definable, even in the bicolour language.
\end{example}
\begin{proof}
Let $\mathbb{B}$ be a de Vries algebra (in the sense of Definition \ref{def1}), containing an element $a$ such that $a \not\prec a $ (that is, by de Vries duality \cite{deVries}, considering a compact Hausdorff space with a regular open set $O$ such that $\overline{O} \not\subseteq O$). 

 Let $I = {\prec}^{ -1}(a)$ and let $\theta$ be its associated congruence (recall the discussion in Definition \ref{defcong}). By axiom (S8), it is clear that $I$ is a round ideal and, hence, that $\theta$ is a white congruence. Moreover, by (S7) and (S8), the filter associated to $I$, i.e. the 1-kernel of $\theta$, is round. It follows that $\theta$ is also a black congruence and, hence, a strong one. Now, we have that $\mathbb{B}/ \theta $ does not satisfy (S5). Indeed, we have $a^\theta \neq 0^\theta$, since otherwise, we would have $a \prec a$. Moreover, if $b^\theta \prec^\theta a^\theta$ for some $b \in \mathbb{B}$, then $b^\theta = 0^\theta$ since 
 
 \begin{center}
 \begin{tabular}{rll}
 &$b^\theta \prec^\theta a^\theta$ \\
 $\Leftrightarrow$ & $ b \prec c \mathrel{\theta} a$& By definition of $\prec^\theta$ (see \eqref{eq_def_prec_theta}) \\
 $\Leftrightarrow$ & $b \prec c$ and $c \vee d = a \vee d$ for some $d \prec a$& \\
 $\Rightarrow$ & $b \prec c$ and $c \vee d = a \vee d$ for some $d \leq a$ & By axiom (S6) \\
 $\Rightarrow$ & $b \prec c$ and $c \leq a$\\
 $\Rightarrow$ & $b \prec a$& By axiom (S4)\\
 $\Leftrightarrow$ & $b^\theta = 0 ^\theta$ & By definition of $\theta$. 
 \end{tabular}
 \end{center}
 \end{proof}
 
 \section*{Conclusions}

 In the present paper, we used previously established dualities (\cite{Quasimod},  \cite{Sourab} and, to a lesser extent, \cite{deVries}) to enthrone subordination algebras as models for tense/modal logics.  Without revoking modal algebras as the suited algebraic formalism, subordination algebras can be seen as parallel models. Indeed, we proved in Section \ref{Section2} two subordination completeness theorems which states, that as far as scheme are concerned, soundness for subordination algebras is equivalent to soundness for modal algebras. With these completeness outcomes in mind, one can therefore enlarge the pool of available counterexamples in the search of a non-provability result. However, this enlarged pool came with a downside: the logic of a class of subordination algebras is not guaranteed to be closed under substitution and hence, is not guaranteed to be normal.  We therefore also provided a functor, from \textbf{Sub} to \textbf{MA}, which can be used to characterise the class of subordination algebras whose logic is normal.
 
Exploring further the relation between subordination algebras and modal logic, we looked in Section \ref{Section_universal algebra} for constructions preserving validity. Since subordination algebras contain some universal algebraic flavour (clearly visible in their alternative presentation as quasi-modal algebras in \cite{Quasimod}), we naturally turned to the notions of subobject, congruence and (Cartesian) product. Some of these notions had been studied by Celani and Castro in \cite{qcong} and \cite{Castro}, and we presented them here in the subordination formalism instead of the quasi-modal one.  These notions leaded us in Section \ref{Section3} to a Birkhoff HSP theorem.

Finally, we established correspondence theorems, via subordination algebras, between three languages: the modal one, the subordination one and the accessibility one. In this paper, we addressed two correspondence theorems: between modal and accessibility first and then between modal and subordination (the third one being discussed for instance in \cite{Balbiani} and \cite{SantoliThese}). The main idea behind the proofs of these theorems is similar to the topological proof of Sahlqvist's theorem in \cite{Sahlqvist}: be in conditions to use an intersection lemma to eliminate variables.  The paper is concluded with some examples of known conditions (which can be found for instance in \cite{Chellas} and \cite{Vakar}) which are untranslatable from one language to another.

\paragraph{Acknowledgements}
We are  thankful to the anonymous referees for their  helpful comments. This has resulted in restructured theorems and proofs, which should be more readable. We also thank Alessandra Palmigiano for her advices on the paper.

%
%



\end{document}